\documentclass[11pt,a4paper,reqno]{amsart}
\usepackage{amsmath,amssymb,amsthm,amscd}
\usepackage[czech,english]{babel}
\usepackage[all]{xy}
\usepackage{hyperref}

\title[On purity and coderived and singularity categories]%
{On purity and applications to coderived and singularity categories}
\date{\today}


\author{Jan \v S\v tov\'\i\v cek}
\address{
Charles University in Prague \\
Faculty of Mathematics and Physics \\
Department of Algebra \\
Sokolovska 83, 186 75 Praha \\
Czech Republic
}
\email{stovicek@karlin.mff.cuni.cz}

\thanks{This research was supported by GA~\v{C}R P201/12/G028.}

\renewcommand{\iff}{if and only if }
\newcommand{\st}{such that }
\newcommand{\wrt}{with respect to }
\newcommand{\wfs}{weak factorization system }
\newcommand{\wfss}{weak factorization systems }
\newcommand{\ep}{\varepsilon}
\newcommand{\ph}{\varphi}

\newcommand{\Z}{\mathbb{Z}}

\newcommand{\m}{\mathfrak{m}}
\DeclareMathOperator{\Hom}{Hom}
\DeclareMathOperator{\HOM}{\mathcal{H}\mathit{om}}

\DeclareMathOperator{\Ext}{Ext}

\DeclareMathOperator{\Ker}{Ker}
\DeclareMathOperator{\Img}{Im}
\DeclareMathOperator{\Coker}{Coker}

\newcommand{\dw}{\mathrm{dw}}
\newcommand{\lifts}{\,\square\,}
\newcommand{\A}{\mathcal{A}}

\newcommand{\C}{\mathcal{C}}

\newcommand{\F}{\mathcal{F}}
\newcommand{\G}{\mathcal{G}}
\newcommand{\I}{\mathcal{I}}
\newcommand{\clL}{\mathcal{L}}
\newcommand{\M}{\mathcal{M}}
\newcommand{\clP}{\mathcal{P}}
\newcommand{\R}{\mathcal{R}}
\newcommand{\clS}{\mathcal{S}}
\newcommand{\T}{\mathcal{T}}
\newcommand{\U}{\mathcal{U}}
\newcommand{\X}{\mathcal{X}}
\newcommand{\Y}{\mathcal{Y}}
\newcommand{\Modu}[1]{{\mathrm{Mod}\textrm{-}{#1}}}
\newcommand{\ModR}{{\mathrm{Mod}\textrm{-}R}}
\newcommand{\modR}{{\mathrm{mod}\textrm{-}R}}
\newcommand{\PProj}[1]{\mathrm{PProj}{#1}}
\newcommand{\Proj}[1]{\hbox{\rm Proj-}{#1}}
\newcommand{\Flat}[1]{\hbox{\rm Flat-}{#1}}
\newcommand{\proj}[1]{\hbox{\rm proj-}{#1}}
\newcommand{\Cot}[1]{\hbox{\rm Cot-}{#1}}
\newcommand{\PInj}[1]{\mathrm{PInj}{#1}}
\newcommand{\Inj}[1]{\mathrm{Inj}\,{#1}}
\newcommand{\dgInj}[1]{\mathrm{dg}\textrm{-}\mathrm{Inj}\,{#1}}
\newcommand{\Fpinj}[1]{{\mathrm{Fpinj}\,{#1}}}
\newcommand{\InjR}{{\mathrm{Inj}R}}
\newcommand{\FpinjR}{{\mathrm{Fpinj}R}}

\newcommand{\Ab}{\mathrm{Ab}}

\newcommand{\Filt}{\mathrm{Filt}}
\newcommand{\Infl}{\operatorname{Infl}}
\newcommand{\Defl}{\operatorname{Defl}}
\newcommand{\fp}{\mathrm{fp}\,}
\newcommand{\Cpx}{\mathbf{C}}
\newcommand{\Cac}{\mathbf{C}_\mathrm{ac}}
\newcommand{\Ccoac}{\mathbf{C}_\mathrm{coac}}
\newcommand{\Cpac}{\mathbf{C}_\mathrm{pac}}

\newcommand{\Htp}{\mathbf{K}}
\newcommand{\Hac}{\mathbf{K}_\mathrm{ac}}

\newcommand{\Hpac}{\mathbf{K}_\mathrm{pac}}

\newcommand{\Der}{\mathbf{D}}
\newcommand{\Dco}{\mathbf{D}^\mathrm{co}}
\newcommand{\Db}{\mathbf{D}^\mathrm{b}}
\newcommand{\Dperf}{\mathbf{D}^\mathrm{perf}}
\newcommand{\Dsing}{\mathbf{D}^\mathrm{sing}}
\newcommand{\Dle}[1]{\mathbf{D}^{\le{#1}}}
\newcommand{\Dge}[1]{\mathbf{D}^{\ge{#1}}}

\newcommand{\Icell}{\I\textrm{-}\mathrm{cell}}
\newcommand{\cof}{\mathrm{Cof}}
\newcommand{\Cof}{\mathcal{C}}

\newcommand{\we}{\mathrm{W}}
\newcommand{\We}{\mathcal{W}}

\newcommand{\fib}{\mathrm{Fib}}
\newcommand{\Fib}{\mathcal{F}}

\newcommand{\Ho}{\operatorname{Ho}}
\newcommand{\la}{\longrightarrow}

\newcommand{\op}{^\mathrm{op}}
\newcommand{\inv}{^{-1}}

\theoremstyle{plain}
\newtheorem{thm}{Theorem}[section]
\newtheorem{prop}[thm]{Proposition}
\newtheorem{lem}[thm]{Lemma}
\newtheorem{cor}[thm]{Corollary}

\newtheorem*{thm_cotorstion-pac}{Theorem~\ref{thm:cotorsion-pac} and Corollary~\ref{lem:recollement-pure}}
\newtheorem*{cor_gor-inj}{Corollary~\ref{cor:gor-inj}}
\newtheorem*{thm_fp-inj-model}{Theorem~\ref{thm:fp-inj-model} and Corollary~\ref{cor:fp-inj}}
\newtheorem*{thm_singularity}{Theorem~\ref{thm:singularity}}

\theoremstyle{definition}
\newtheorem{defn}[thm]{Definition}
\newtheorem{hyp}[thm]{Hypothesis}

\theoremstyle{remark}
\newtheorem{rem}[thm]{Remark}
\newtheorem{expl}[thm]{Example}
\newtheorem{warn}[thm]{Warning}

\begin{document}

\begin{abstract}
Given a locally coherent Grothendieck category $\G$, we prove that the homotopy category of complexes of injective objects (also known as the coderived category of $\G$) is compactly generated triangulated. Moreover, the full subcategory of compact objects is none other than $\Db(\fp\G)$. If $\G$ admits a generating set of finitely presentable objects of finite projective dimension, then also the derived category of $\G$ is compactly generated and Krause's recollement exists. Our main tools are (a) model theoretic techniques and (b) a systematic study of the pure derived category of an additive finitely accessible category.
\end{abstract}

\maketitle

\setcounter{tocdepth}{1}
\tableofcontents

\section{Introduction}
\label{sec:intro}

The aim of this paper is to establish (at first perhaps surprising) finiteness results for coderived and singularity categories of locally coherent Groth\-end\-ieck categories. Similar theory for locally noetherian Grothendieck categories have been developed in the last decade by several authors, motivated by the Grothendieck duality and Gorenstein homological algebra. When one deals with functor categories in representation theory or the theory of triangulated categories, one, however, often encounters a non-noetherian categories. Treating the locally coherent situation is considerably more involved than the locally noetherian one and requires a systematic use of purity and model theoretic techniques.

\smallskip

To start with, it was known since the dawn of triangulated categories in algebra in the 1960's that computations in bounded derived categories could be done in homotopy categories of complexes. One just needed to replace each complex by an injective resolution. This is at least theoretically much easier since one avoids the non-commutative localization otherwise inherent to the construction of derived categories. The situation with unbounded derived categories was traditionally much more puzzling. It was independently resolved by Spaltenstein~\cite{Spa88} and Alonso Tarr\'{\i}o, Jerem\'{\i}as L\'{o}pez and Souto Salorio~\cite{AJS00} in the context of algebraic geometry, and by Joyal~\cite{Joy84} and Beke~\cite{Be00} (see also~\cite[Example 3.2]{Hov02-models-cot}) with homotopy theory in mind.
The problem is that not every unbounded complex of injectives qualifies as a resolution---there are in general acyclic complexes of injectives which are not contractible. Nevertheless, it was proved that given a Grothendieck category $\G$, the unbounded derived category $\Der(\G)$ was equivalent to a full subcategory of the homotopy category $\Htp(\Inj\G)$. The conclusion by then was that $\Htp(\Inj\G)$ was not quite the correct category for homological algebra and this full subcategory should have been considered instead.

This point of view has radically changed with the work of Krause \cite{Kr05} and J\o{}rgensen \cite{Jo05}. They gave clear interpretations to the homotopy categories of complexes of injective sheaves and of projective modules. We refer to~\cite{Orl04,Kr05,Nee10-modern} for details and applications in representation theory, algebraic geometry and theoretical physics. The present work can be viewed as a continuation of this effort, that is of \cite{Kr05,Jo05} as well as Neeman's \cite{Nee08,Nee10,Nee14} and Becker's \cite{Beck14}. The main emphasis is put on \emph{non-noetherian} situations.

\smallskip

Let $\G$ be a locally finitely presentable Grothendieck category. This includes for instance
\begin{enumerate}
\item[(a)] module categories over rings or small preadditive categories,
\item[(b)] nice localizations of module categories, relevant in particular in model theory~\cite{CB94,Her97,Kr97,Prest09} and the theory of triangulated categories~\cite{Nee97,Nee01,ChKeNee01,Bel00,DASS14}, and
\item[(c)] various categories of quasi-coherent sheaves~\cite{Th90}.
\end{enumerate}
Such categories arising in case (b) are often \emph{locally coherent}, i.e.\ the full subcategory $\fp\G$ of finitely presentable objects is abelian, but not necessarily noetherian. It still makes perfect sense to study the bounded derived category $\Db(\fp\G)$. Although it is equivalent to a full subcategory of $\Der(\G)$, it is well known that the relation between $\Db(\fp\G)$ and $\Der(\G)$ is not the same as the one between $\fp\G$ and $\G$.
Namely, $\Der(\G)$ is frequently a compactly generated triangulated category (an analog of locally finitely presentable in the triangulated world), generated by a set of objects from $\Db(\fp\G)$. At the same time, however, often not every object of $\Db(\fp\G)$ is compact in $\Der(\G)$.

Here is the main contribution of~\cite{Kr05}. He proved that if $\G$ is locally noetherian, then $\Db(\fp\G)$ can be identified with the compact objects of $\Htp(\Inj\G)$. Out first main result extends this to locally coherent Grothendieck categories. In fact, we also show an example of a non-coherent ring where this fails, so we cannot just drop the coherence assumption.

\begin{thm_fp-inj-model}
Given a locally coherent Grothendieck category $\G$, then $\Htp(\Inj\G)$ is a compactly generated triangulated category. Moreover, a functor
\[ I\colon \Db(\fp\G) \la \Htp(\Inj\G) \]
which assigns to each bounded complex $X \in \Db(\fp\G)$ a fixed injective resolution $I(X)$ of $X$ and similarly for morphisms, restricts to an equivalence between $\Db(\fp\G)$ and the full subcategory $\Htp(\Inj\G)^c$ of compact objects in $\Htp(\Inj\G)$.
\end{thm_fp-inj-model}

For a general Grothendieck category, $\Htp(\Inj\G)$ is only well generated in the sense of \cite{Nee01} by~\cite[Theorem 3.13]{Nee14}. Our result is a tremendous refinement of this---it allows us to go back from the world of $\mu$-compact generators for a regular cardinal $\mu$ to classical compact objects. It also resolves Krause's \cite[Conjecture 1]{Kr12}.

\smallskip

There is usually a twofold motivation for studying the homotopy category of injective objects. Firstly, it provides a novel point of view at the Grothendieck duality in algebraic geometry; see~\cite{KrIy06,Nee08}. Secondly, many problems in representation theory, especially when one encounters various flavors of Gorenstein homological algebra or Tate cohomology as in~\cite{Bu86,Rick89,Orl04,Kr05}, lead to cosyzygies of infinite order. That is, we encounter objects $X \in \G$ which admit an exact resolution
\[ \cdots \la I_2 \la I_1 \la I_0 \la X \la 0 \]
where all $I_n$, $n \ge 0$, are \emph{injective}. This for instance also happens in~\cite{DASS14}, a work in progress and a direct motivation for this research, aiming at explaining refined versions of Brown-Adams representability in Kasparov's bivariant $K$-theory of $C^\ast$-algebras~\cite{Koeh11,MeNe12}. Somewhat surprisingly, we show that there is a qualitative difference between $N$-th cosyzygies for $N$ large but finite, and cosyzygies of infinite order (see also Remark~\ref{rem:vNreg}). This seems to be a new result even for noetherian rings.

\begin{cor_gor-inj}
An infinite order cosyzygy module is always cotorsion.
\end{cor_gor-inj}

The class of cotorsion modules can be viewed as the smallest extension of the class pure-injective modules which can be defined in terms of usual homological functors; see~\cite[\S4.6]{Prest09}. Although somewhat technical in nature, this result is extremely useful when dealing with homological algebra for infinitely generated modules. To get an idea what kind of problems one encounters with non-cotorsion modules, we refer to~\cite{EST04}.

A noteworthy consequence of Corollary~\ref{cor:gor-inj} is that it allows to work with non-noetherian analogs of Gorenstein rings exactly as Buchweitz did in~\cite{Bu86} with noetherian rings. For instance, we shall prove in Proposition~\ref{prop:ding-chen} that Gillespie's Ding-Chen injective modules~\cite{Gil10} are precisely cosyzygies of infinite order. Another application of similar flavor will be given in~\cite{DASS14}.

As is known from Buchweitz' work~\cite{Bu86}, the stable module categories of Iwanaga-Gorenstein rings are equivalent to what is nowadays called singularity categories~\cite{Orl04}. These are usually defined as the Verdier quotients $\Dsing(\fp\G) = \Db(\fp\G)/\Dperf(\fp\G)$, where $\Dperf(\fp\G)$ is the class of perfect complexes, i.e.\ the compact objects of $\Der(\G)$. As in the locally noetherian case~\cite{Kr05}, we obtain a recollement involving the `inductive completion' of $\Dsing(\fp\G)$ also in the locally coherent case.

\begin{thm_singularity}
Let $\G$ be a locally coherent Grothendieck category which admits a generating set of finitely presentable objects of finite projective dimension. Then there is a recollement
\[
\xymatrix@1{
\Hac(\Inj\G) \ar[rr]|(.45)\hole|{i_*} &&
\Htp(\Inj\G) \ar@/_3ex/[ll]_{i^*} \ar@/^3ex/[ll]^{i^!} \ar[rr]|(.45)\hole|{Q} &&
\Der(\G) \ar@/_3ex/[ll]_{L} \ar@/^3ex/[ll]^{R}
},
\]
and both $\Der(\G)$ and $\Hac(\Inj\G)$ are compactly generated. Moreover, the full subcategory of compact objects of $\Hac(\Inj\G)$ is equivalent to the idempotent completion of the singularity category $\Dsing(\fp\G)$.
\end{thm_singularity}

Note that for coherent rings, the theorem in particular guarantees the existence of non-trivial cosyzygies of infinite order provided that the global dimension of $\modR$, the category of finitely presented modules, is infinite.

\smallskip

Let us now shortly mention the techniques used to prove our results. Although the results look formally identical to those in~\cite{Kr05}, the generalization from the locally noetherian to the locally coherent situation makes the proofs considerably more involved. We systematically use and combine
\begin{enumerate}
\item model theoretic techniques \cite{Hov02-models-cot,Gil04,Beck14} and,
\item purity and the related pure homological algebra.
\end{enumerate}

The rationale behind (1) is that several predecessors of the present paper, including~\cite{BEIJR12,Mur,Nee08,Nee10,Nee14}, use cotorsion pairs and the related approximation theory. In view of \cite{Hov02-models-cot}, this is none other than implicitly using certain abelian model structures. Doing so systematically, many arguments are much more transparent and conceptual. For instance, we often follow Positselski~\cite{Pos11} and Becker~\cite{Beck14}, and view the homotopy category $\Htp(\Inj\G)$ as the localization of $\Cpx(\G)$, where we formally invert a certain precisely specified subclass of the class of quasi-isomorphisms. Of course, the real strength of the approach is unleashed when one combines the model theoretic techniques with a wide repertoire of algebraic tricks for computing with complexes.

The second main component consists of a deeper look at the pure derived category, which has been formally considered for instance in \cite[\S5.3]{ChHo02} and \cite{Kr12}. Here we can express our main contribution in terms of another recollement:

\begin{thm_cotorstion-pac}
Let $\A$ be an additive finitely accessible category with the standard pure exact structure. Then the unbounded derived category $\Der(\A)$ admits a recollement
\[
\xymatrix@1{
\Hpac(\A) \ar[rr]|(.45)\hole|{i_*} &&
\Htp(\A) \ar@/_3ex/[ll]_{i^*} \ar@/^3ex/[ll]^{i^!} \ar[rr]|(.45)\hole|{j^*} &&
\Der(\A) \ar@/_3ex/[ll]_{j_!} \ar@/^3ex/[ll]^{j_*}
}.
\]
Moreover, the essential image of $j_!$ is $\Htp(\PProj\A)$, the homotopy category of pure projective objects, and the essential image of $j_*$ is $\Htp(\PInj\A)$, the homotopy category of pure injective objects.
\end{thm_cotorstion-pac}

As an \emph{additive finitely accessible category} is merely a synonym for the category of flat modules over a small additive category, the existence of the recollement and the description of the essential image of $j_!$ has been already obtained by Neeman in~\cite[Theorem 8.6]{Nee08} and~\cite[Theorem 0.1]{Nee10}. What is new here and crucial in our situation is the description of the essential image of $j_*$.

We conclude the introduction by mentioning consequences of this last result. Given a ring $R$ and $\A = \Flat{R}$, then $j_*j^*(R)$ is just a (minimal) proper cotorsion resolution of $R$. If $R$ is left coherent, it is an immediate consequence of~\cite[Lemma 3.1.6]{Xu96} that flat cotorsion modules are pure injective. In particular, $j_*j^*(R)$ is also a (minimal) pure injective resolution of $R$ (see also~\cite{En84}). This answers the question posed in the introductions to \cite{Nee08,Nee10} on what $j_*j^*(R)$ looks like, motivated by the observation that if $R$ is commutative with a dualizing complex $I$, then $j_*j^*(R)$ is also homotopy isomorphic to $\HOM_R(I,I)$.

Secondly, the theorem can be interpreted as that there is nothing like a pure singularity category or pure Gorenstein homological algebra. Indeed, any pure acyclic complex of pure injective objects is necessarily contractible. Interestingly, this holds despite the fact that the pure global dimension of $\A$ may well be infinite. Indeed, there exists a valuation domain $R$ whose maximal ideal $\m$ has infinite projective dimension by~\cite[Proposition 2.62]{Os73}, and so $\m$ has infinite pure projective dimension in $\A = \Flat{R}$.

\section{Weak factorization systems and cotorsion pairs}
\label{sec:wfss-cot}

In this section we recall some facts about the primitives of Quillen model structures---\wfss and complete cotorsion pairs. We shall use these concepts frequently in the remaining text and, although some of the abstract results can be stated more generally, we will always work in the framework of accessible categories. To this end, let $\kappa$ be a regular cardinal and $\A$ be a category with arbitrary $\kappa$-direct limits, i.e.\ colimits of all diagrams whose shapes are $\kappa$-directed posets. We shall use the convention of~\cite{AR94,GU71} where a poset $I$ is $\kappa$-directed if each subset of $I$ of cardinality \emph{smaller} than $\kappa$ has an upper bound. Recall that an object $F \in \A$ is called \emph{$\kappa$-presentable} if for each $\kappa$-direct system $(Y_i \mid i \in I)$ in $\A$, the canonical map
\[ \varinjlim \A(F,Y_i) \la \A(F,\varinjlim Y_i) \]
is an isomorphism.

\begin{defn}[\cite{AR94,GU71}] \label{defn:fin-acc}
If $\A$ has $\kappa$-direct limits and there exists a small set $\clS$ of $\kappa$-presentable objects \st every $X \in \A$ can be expressed as a $\kappa$-direct limit of objects from $\clS$, then $\A$ is called \emph{$\kappa$-accessible}. If, moreover, $\A$ is cocomplete, it is called \emph{locally $\kappa$-presentable}.

In the special (and prominent) case where $\kappa = \aleph_0$, we speak of \emph{finitely presentable} objects, and of \emph{finitely accessible} or \emph{locally finitely presentable} categories. The full subcategory of finitely presentable objects of $\A$ is then denoted by $\fp\A$.
\end{defn}

\subsection{Weak factorization systems}
\label{subsec:wfss}

Let us remind the reader of the lifting relation between morphisms. Given morphisms $f\colon A \to B$ and $g\colon X \to Y$ in $\A$. we say that $f$ has the \emph{left lifting property} for $g$, or that $g$ has the \emph{right lifting property} for $f$, if for any commutative square given by the solid arrows
\[
\xymatrix{
A \ar[r] \ar[d]_f     & X\phantom{,} \ar[d]^g      \\
B \ar[r] \ar@{.>}[ur] & Y,
}
\]
there exists a diagonal dotted morphism \st both triangles commute. In symbols we write we write $f \lifts g$.
Note that if $\A$ is additive and we view the maps $f$ and $g$ as two-term complexes, $f \lifts g$ simply means that any chain map from $f$ to $g$ is null-homotopic.

\begin{defn} \label{defn:wfs}
Let $(\clL,\R)$ be a pair of classes of morphisms in $\A$. We say that $(\clL,\R)$ is a \emph{\wfs}if
\begin{enumerate}
\item $\clL$ and $\R$ are closed under retracts.
\item $f \lifts g$ for all $f \in \clL$ and $g \in \R$.
\item For every morphism $h\colon X \to Y$ in $\A$, there is a factorization
\[
\xymatrix@R=1em{
X \ar[rr]^h \ar@/_/[dr]_f && Y  \\
& Z \ar@/_/[ur]_g
}
\]
with $f \in \clL$ and $g \in \R$.
\end{enumerate}

We say that $(\clL,\R)$ is a \emph{functorial \wfs}if the factorization in~(3) can be chosen functorially in $h$.
\end{defn}

Beware that the factorizations as in (3), be they functorial or not, need not be unique. This is the reason for the adjective \emph{weak}. The classes $\clL$ and $\R$ of a \wfs have well-known closure properties which we recall next, along with the necessary definition of a transfinite composition of morphisms:

\begin{defn}[\cite{Hir03,Hov99}] \label{defn:lambda-seq}
Let $\lambda$ be an ordinal number and $(X_\alpha, f_{\alpha\beta} \mid \alpha<\beta<\lambda)$ be a direct system in $\A$ indexed by the ordinal numbers below~$\lambda$. Such a system is called a \emph{$\lambda$-sequence} provided that for each limit ordinal $\mu<\lambda$, the subsystem $(X_\alpha, f_{\alpha\beta} \mid \alpha<\beta<\mu)$ has colimit $X_\mu$, the colimit morphisms being $f_{\alpha\mu}\colon X_\alpha \to X_\mu$. In other words, we require that $X_\mu \cong \varinjlim_{\alpha<\mu} X_\alpha$ canonically.

The \emph{composition} of a $\lambda$-sequence $(X_\alpha, f_{\alpha\beta} \mid \alpha<\beta<\lambda)$ is defined as the colimit morphism $X_0 \to \varinjlim_{\alpha<\lambda} X_\alpha$, if the colimit exists.
Finally given a class $\I$ of morphisms in $\A$, then the \emph{transfinite compositions} of morphisms of $\I$ are by definition the compositions of those $\lambda$-sequences $(X_\alpha, f_{\alpha\beta} \mid \alpha<\beta<\lambda)$ with $f_{\alpha,\alpha+1} \in \I$ for each $\alpha+1<\lambda$.
\end{defn}

\begin{lem} \label{lem:wfss-closure}
Let $(\clL,\R)$ be a \wfs in $\A$. Then $\clL = \{ f \mid f \lifts g \textrm{ for all } g \in \R \}$ and $\R   = \{ g \mid f \lifts g \textrm{ for all } f \in \clL \}$. Moreover,
\begin{enumerate}
\item $\clL$ is closed under pushouts. That is, a pushout of any $f \in \clL$ along an arbitrary morphism of $\A$ again belongs to $\clL$. Dually, $\R$ is closed under pullbacks.
\item $\clL$ is closed under transfinite compositions. More precisely, any transfinite composition of morphisms of $\clL$ again belongs to $\clL$.
\end{enumerate}
\end{lem}

\begin{proof}
This is easy and well known.
\end{proof}

The lemma motivates the following definition, where we apply all the above mentioned closure properties of $\clL$ to a given class of morphisms $\I$.

\begin{defn}[\cite{Hir03,Hov99}] \label{defn:I-cell}
If $\I$ is a class of morphisms of $\A$, then a \emph{relative $\I$-cell complex} is a transfinite composition of pushouts of morphisms from $\I$. We denote the class of relative $\I$-cell complexes by $\Icell$.
\end{defn}

Finally, there is an important technical result, based on so-called Quillen's small object argument, which explains how \wfss are generated by small sets of morphisms from the left. In particular it shows that under mild (albeit somewhat technical) assumptions, we have identified \emph{all} relevant closure properties of $\clL$ in Lemma~\ref{lem:wfss-closure}.

\begin{prop} \label{prop:Quillen-soa}
Let $\A$ be a $\kappa$-accessible category for some $\kappa$, and let $\M$ be a class of morphisms \st
\begin{enumerate}
\item arbitrary pushouts of morphisms of $\M$ exist and belong again to $\M$,
\item arbitrary coproducts of morphisms of $\M$ exist and belong to~$\M$,
\item arbitrary transfinite compositions of maps of $\M$ exist and are in~$\M$.
\end{enumerate}
Let $\I \subseteq \M$ be a small set of morphisms and put
\[
\R   = \{ g \mid f \lifts g \textrm{ for all } f \in \I \}
\quad\textrm{and}\quad
\clL = \{ f \mid f \lifts g \textrm{ for all } g \in \R \}.
\]
Then $(\clL,\R)$ is a functorial \wfs in $\A$ and $\clL$ consists precisely of retracts of relative $\I$-cell complexes.
\end{prop}

\begin{proof}
We refer to~\cite[Proposition 10.5.16 and Corollary 10.5.22]{Hir03} or to \cite[Theorem 2.1.14 and Corollary 2.1.15]{Hov99}. Although the assumptions in~\cite{Hir03,Hov99} differ slightly, identical proof works.
\end{proof}

\subsection{Cotorsion pairs in exact categories}
\label{subsec:cotorsion}

If our category $\A$ carries more structure, say it is an exact or abelian category, we can often work with objects of $\A$ rather than morphisms. This is often more convenient. The lifting relation between morphisms from \S\ref{subsec:wfss} is then replaced by the orthogonality with respect to $\Ext^1$, and \wfss transform to so called complete cotorsion pairs. This reduction, which we are going to briefly explain, was worked out by Hovey~\cite{Hov02-models-cot} while similar considerations in a special case appeared at around the same time in~\cite{Ros02}. For the further developments we refer to~\cite{St14} and the references therein.

Let us fix some notation first. Given a skeletally small preadditive category $\clS$, by a \emph{right $\clS$-module} we mean an additive functor $\clS\op \to \Ab$. We denote the category of all right $\clS$-modules by $\Modu\clS$. If $\clS$ has only one object with endomorphism ring $R$, then $\Modu\clS$ is the same as $\ModR$. Module categories are prototypical examples where the reduction from \wfss to complete cotorsion pairs works.

We will, however, need to work in a more general setup.
Recall from~\cite{Bu10} that an \emph{exact category} $\A$ is an additive category with a designated collection of kernel-cokernel pairs $(i,d)$ sharing certain properties with short exact sequences in abelian categories. Given such a kernel-cokernel pair $(i,d)$ with $i\colon K \to L$ and $d\colon L \to M$, we will usually write it in the form 
\[ 0 \la K \stackrel{i}\la L \stackrel{d}\la M \la 0 \]
and, following~\cite[Appendix A]{Kel90}, we call it a \emph{conflation}. The map $i$ is then called an \emph{inflation} and the map $d$ a \emph{deflation}. 

If $\A$ is an exact category, we can define Yoneda Ext functors $\Ext^i_\A\colon \A\op \times \A \to \Ab$ using the conflations as if they were short exact sequences. To be precise, the only obstacle which would prevent us from doing so is that $\Ext^i_\A(X,Y)$ may not be a small set in general, but this pitfall fortunately never occurs for our choices of $\A$.
It is well known that these Ext functors have the usual properties which one would expect; see for instance~\cite[\S5]{St14}. If $\A = \Modu\clS$ or $\ModR$ with the \emph{abelian exact structure} (i.e.\ conflations are chosen to be precisely the short exact sequences), we will simply write $\Ext^i_\clS$ or $\Ext^i_R$, respectively. Now we can give a key definition:

\begin{defn} \label{defn:cotorsion}
Let $\A$ be an exact category. A pair $(\X,\Y)$ of full subcategories of $\A$ is called a \emph{cotorsion pair} provided that
\begin{align*}
\X &= \{X \in \A \mid \Ext^1_\A(X,Y) = 0 \textrm{ for all } Y \in \Y\} \quad \textrm{and} \\
\Y &= \{Y \in \A \mid \Ext^1_\A(X,Y) = 0 \textrm{ for all } X \in \X\}.
\end{align*}

A cotorsion pair $(\X,\Y)$ is said to be \emph{complete} if every $A \in \A$ admits so-called \emph{approximation sequences}; that is, conflations of the form
\[
\ep_A\colon \; 0 \to A \to Y_A \to X_A \to 0 \quad \textrm{ and } \quad
\ep^A\colon \; 0 \to Y^A \to X^A \to A \to 0
\]
with $X_A, X^A \in \X$ and $Y_A, Y^A \in \Y$.

The cotorsion pair is \emph{functorially complete} if the approximation sequences $\ep_A$ and $\ep^A$ can be taken functorial in $A$.
\end{defn}

We again stress that the approximation sequences, even if they are functorial, are usually very non-unique. The relation to \wfss is very straightforward.

\begin{prop} \label{prop:cot-to-wfs}
Let $(\X,\Y)$ be a complete cotorsion pair in an exact category $\A$. Denote by $\Infl\X$ the class of all inflations with cokernels in $\X$, and by $\Defl\Y$ the class of all deflations with kernels in $\Y$. Then $(\Infl\X,\Defl\X)$ is a \wfs in $\A$, and it is functorial \iff $(\X,\Y)$ is such.
\end{prop}

\begin{proof}
This result was essentially obtained by Hovey in \cite[\S5]{Hov02-models-cot}. An explicit statement of this form is given in~\cite[Theorem 5.13]{St14}.
\end{proof}

In other words, Proposition~\ref{prop:cot-to-wfs} says that $(\X,\Y) \mapsto (\Infl\X,\Defl\Y)$ is an injective mapping from the class of (functorially) complete cotorsion pairs to the class of (functorial) \wfss in $\A$. It is not difficult to describe the image of this mapping (see again \cite{Hov02-models-cot} or \cite{St14}), but we shall only briefly mention this in Section~\ref{sec:models}. Here we shall rather focus on the analogs of Lemma~\ref{lem:wfss-closure} and Proposition~\ref{prop:Quillen-soa} for cotorsion pairs.

As for Lemma~\ref{lem:wfss-closure}, only part~(2) is relevant here. It leads to a useful formalization of transfinite extensions in exact categories which, following~\cite{GT}, we call filtrations.

\begin{defn} \label{defn:filtr}
Let $\A$ be an exact category and $\clS$ be a class of objects in~$\A$. An \emph{$\clS$-filtration} is a $(\lambda+1)$-sequence $(X_\alpha, i_{\alpha\beta} \mid \alpha<\beta\le\lambda)$ indexed by an ordinal successor $\lambda+1$ and satisfying the following three conditions:
\begin{enumerate}
\item $X_0 = 0$,
\item $i_{\alpha\beta}\colon X_\alpha \to X_\beta$ is an inflation in $\A$ for each $\alpha<\beta\le\lambda$,
\item $\Coker i_{\alpha,\alpha+1} \in \clS$ up to isomorphism for each $\alpha<\lambda$.
\end{enumerate}
An object $X \in \A$ is \emph{$\clS$-filtered} if there is an $\clS$-filtration $(X_\alpha, i_{\alpha\beta} \mid \alpha<\beta\le\lambda)$ with $X = X_\lambda$. The class of all $\clS$-filtered objects in $\A$ is denoted by $\Filt(\clS)$.

A class $\C \subseteq \A$ is called \emph{deconstructible} if it is of the form $\C = \Filt(\clS)$ for a \emph{small} set $\clS$ of objects of $\A$.
\end{defn}

Now we can state the analog of Lemma~\ref{lem:wfss-closure}(2) which is known as the Eklof Lemma (see~\cite[Lemma 3.1.2]{GT}) and which essentially says that the left hand side class of any cotorsion pair is closed under transfinite extensions.

\begin{prop}[Eklof Lemma] \label{prop:eklof}
Let $\A$ be an exact category, $\clS$ be a class of objects and $Y \in \A$ be \st $\Ext^1_\A(S,Y) = 0$ for all $S \in \clS$. If $X$ is $\clS$-filtered, then $\Ext^1_\clS(X,Y) = 0$.
\end{prop}

\begin{proof}
This is a little more general variant of~\cite[Proposition 2.12]{SaSt11}. Let $X,Y,\clS$ be as in the statement and let $(X_\alpha, i_{\alpha\beta} \mid \alpha<\beta\le\lambda)$ be an $\clS$-filtration for $X$. Given a conflation
\[ \ep\colon\quad 0 \la Y \stackrel{j}\la E \la X \la 0 \]
in $\A$, our task is to prove that $\ep$ splits. First of all we claim that there is a direct system $(E_\alpha, j_{\alpha\beta} \mid \alpha<\beta\le\lambda)$ which satisfies conditions (2) and (3) of Definition~\ref{defn:filtr} and \st $j_{0\lambda} = j\colon Y \to E$. As expected, we construct such a system by taking the pullbacks of $\ep$ along $i_{\alpha\lambda}$. In the exact categories $\A$ which we study in this paper, such as module categories, this is a nearly trivial claim, but we also give a (rather tedious) abstract proof in Proposition~\ref{prop:rel_filtr} at the level of generality which is stated here.

Now we use transfinite induction on $\alpha$ to construct a cocone $(\ph_\alpha\colon E_\alpha \to Y \mid \alpha\le\lambda)$ of $(E_\alpha, j_{\alpha\beta} \mid \alpha<\beta\le\lambda)$ \st $\ph_0 = 1_Y$. At limit steps we simply take the colimit maps of $(\ph_\beta \mid \beta<\alpha)$. At successor steps, we use the fact that $\ph_\alpha$ extends along $j_{\alpha,\alpha+1}$ since $\Ext^1_\A(\Coker j_{\alpha,\alpha+1},Y) = 0$ by assumption. Finally one sees that $\ph_\lambda\colon E \to Y$ is a retraction of $j\colon Y \to E$.
\end{proof}

We conclude the section with an analog of Proposition~\ref{prop:Quillen-soa} for cotorsion pairs. It also involves somewhat technical assumptions, but these are broadly satisfied. As before, crucial is the last sentence of the statement which describes all the objects in the left hand side class of a cotorsion pair generated by a small set.

\begin{prop} \label{prop:complete-cotor-A}
Let $\A$ be an exact category \st
\begin{enumerate}
\item $\A$ is $\kappa$-accessible for some $\kappa$, and
\item arbitrary transfinite compositions of inflations exist and are again inflations.
\end{enumerate}
Suppose further that $\clS \subseteq \A$ is a small set of objects which is generating. That is, we require that for each $X \in \A$ there exist a deflation $\coprod_{i\in I} S_i \to X$ with $S_i \in \clS$ for all $i \in I$. Then
\begin{align*}
\Y &= \{ Y \in \A \mid \Ext^1_\A(S,Y) = 0 \textrm{ for all } S \in \clS \} \qquad \textrm{and} \\
\X &= \{ X \in \A \mid \Ext^1_\A(X,Y) = 0 \textrm{ for all } Y \in \Y   \} 
\end{align*}
form a functorially complete cotorsion pair $(\X,\Y)$ in $\A$. Moreover, $\X$ consists precisely of summands of $\clS$-filtered objects.
\end{prop}

\begin{proof}
This is an instance of \cite[Corollary 2.15]{SaSt11} or \cite[Theorem 5.16]{St14}, so we only need to check the assumptions of those results. In other words, we must check that $\A$ is an efficient exact category in the sense of~\cite[Definition 2.6]{SaSt11}.

The latter amounts to checking axioms (Ax1)--(Ax3) of \cite[Definition 2.6]{SaSt11}, where (Ax1) is precisely our assumption (2). Further, (Ax2) requires that each object of $\A$ be small with respect to the class of inflations (we refer to \cite{SaSt11} for an explanation of this technical condition). This is a consequence of abstract theory of accessible categories, see~\cite[Proposition 1.16 and Remark 2.2(3)]{AR94}, and as such follows from our assumption (1). Finally, proving (Ax3) reduces in view of \cite[Proposition 2.7]{SaSt11} to proving that
\begin{enumerate}
 \item[(a)] $\A$ has a small set of generators,
 \item[(b)] $\A$ has arbitrary coproducts, and
 \item[(c)] every section in $\A$ has a cokernel, i.e.~$\A$ is weakly idempotent complete in the sense of~\cite{Bu10}.
\end{enumerate}

Here (a) is clearly satisfied as we assume that we have a generating small set $\clS$, and (b) follows by \cite[Lemma 1.4]{SaSt11} from our assumption (2). As for (c), note that each $\kappa$-accessible category is even idempotent complete by \cite[Observation 2.4 and Remark 1.21]{AR94}.
\end{proof}

\section{Exact and abelian model structures}
\label{sec:models}

Now we briefly recall the formalism of exact model categories. This is a very natural framework for arguments involving transfinite homotopy colimits which cannot be directly performed at the level of triangulated categories.

Model structures and model categories~\cite{Hir03,Hov99} have been discovered to deal with localizations of categories. If we start with a category $\A$ and a class $\we$ of morphisms in $\A$, we can in principle always construct (up to a possible set theoretic difficulty) a universal functor $Q\colon \A \to \A[\we\inv]$ which makes the maps in $\we$ invertible; see~\cite[\S1.1]{GZ67}. A much harder task is, of course, to understand what $\A[\we\inv]$ looks like. The classical topologically motivated idea behind Quillen's model structures says that if we manage to find more structure, namely two additional classes of morphisms $\cof$ and $\fib$ satisfying the axioms below, the situation becomes much easier.

To fix basic notation and terminology, we define a \emph{model structure} on $\A$ is a triple of classes of morphisms $(\cof,\we,\fib)$ \st
\begin{enumerate}
\item $(\cof,\we\cap\fib)$ and $(\cof\cap\we,\fib)$ are functorial weak factorization systems, and
\item $\we$ is closed under retracts and satisfies the 2-out-of-3 property for compositions.
\end{enumerate}
One easily checks that this is equivalent to the standard definitions in~\cite[\S1.1]{Hov99} or~\cite[\S7.1]{Hir03}.
The morphisms in $\cof,\we,\fib$ are called \emph{cofibrations}, \emph{weak equivalences} and \emph{fibrations}, respectively. If $\A$ is a pointed category, then an object $X \in \A$ is called
\begin{itemize}
 \item \emph{cofibrant} if $0 \to X$ is a cofibration;
 \item \emph{fibrant} if $X \to 0$ is a fibration; and
 \item \emph{trivial} if $0 \to X$ (or equivalently $X \to 0$) is a weak equivalence.
\end{itemize}

If $\A$ is an exact category, one can often encode a model structure in three classes classes of objects rather than morphisms. This is again due to Hovey and the principle is exactly the same as with \wfss and complete cotorsion pairs in Section~\ref{sec:wfss-cot}. For this to work, the model structure has to satisfy a straightforward compatibility condition with the exact structure.

\begin{defn}[\cite{Hov02-models-cot}] \label{defn:exact-model}
Let $\A$ be an exact category. A model structure $(\cof,\we,\fib)$ on $\A$ is \emph{exact} if cofibrations are precisely inflations with cofibrant cokernel and fibrations are precisely deflations with fibrant kernel. If $\A$ is in fact an abelian category considered with the abelian exact structure, we call also the model structure \emph{abelian}.
\end{defn}

In the rest of the text we will exclusively consider exact and abelian model structures. In order to state the description of exact model structures by cotorsion pairs, we need to restrict to weakly idempotent complete exact categories. We remind that this means that every section in $\A$ has a cokernel, or equivalently that every retraction in $\A$ has a kernel (see~\cite[\S7]{Bu10}).

\begin{prop} \label{prop:exact-model}
If $\A$ is a weakly idempotent complete exact category, there is a bijective correspondence between exact model structures on $\A$ and the triples of classes of objects $(\Cof,\We,\Fib)$ satisfying the following conditions:
\begin{enumerate}
\item $(\Cof,\We\cap\Fib)$ and $(\Cof\cap\We,\Fib)$ are functorially complete cotorsion pairs in $\A$.
\item $\We$ is closed under retracts and satisfies the 2-out-of-3 property for extensions. That is, if $0 \to X \to Y \to Z \to 0$ is a conflation in $\A$ and two of $X,Y,Z$ are in $\We$, then so is the third.
\end{enumerate}
To $(\Cof,\We,\Fib)$ one assigns the unique exact model structure $(\cof,\we,\fib)$ on $\A$ whose class of cofibrant objects is $\C$ and whose class of fibrant objects is $\F$. Moreover, a morphism $w$ in $\A$ becomes a weak equivalence \iff $w = w_dw_i$ for an inflation $w_i$ with cokernel in $\We$ and a deflation $w_d$ with kernel in $\We$.
\end{prop}

\begin{proof}
This is \cite[Theorem 3.3]{Gil11}. A detailed account on the same result is also given in~\cite[\S6.2]{St14}. The main idea, however, comes already from~\cite{Hov02-models-cot}.
\end{proof}

Now we shall discuss a few computational aspects of the \emph{homotopy category} $\Ho\A$ of $\A$, which is the traditional name for the localized category $\A[\we\inv]$.

\begin{prop} \label{prop:hom_in_htpy}
Suppose that $\A$ is a weakly idempotent complete exact category and $(\Cof,\We,\Fib)$ represents an exact model structure on $\A$. Then $\Ho\A$ is an additive category and the localization functor $Q\colon \A \to \Ho\A$ is an additive functor. Moreover, whenever $X \in \A$ is cofibrant and $Y \in \A$ is fibrant, $Q$ induces a surjective group homomorphism
\[ \A(X,Y) \stackrel{Q}\la \Ho\A(QX,QY) \]
whose kernel consists precisely of those morphisms $f\colon X \to Y$ which factor through an object in $\Cof\cap\We\cap\Fib$.
\end{prop}

\begin{proof}
The fact that $Q$ induces a surjective map whenever $X$ is cofibrant and $Y$ is fibrant is classical; see~\cite[Theorem 1.2.10]{Hov99}. In fact, \cite[Theorem 1.2.10]{Hov99} uses slightly stronger assumptions on~$\A$ which are unnecessary in our case as explained in~\cite[Remark 6.2]{St14}. The description of when precisely $Qf = Qg$ for $f,g\colon X \to Y$ is obtained in~\cite[Proposition 4.4]{Gil11}. In particular, the preimages of singletons of $\Ho\A(QX,QY)$ are cosets of a fixed subgroup of $\A(X,Y)$, so we can equip $\Ho\A(QX,QY)$ with a unique structure of an abelian group making $\A(X,Y) \to \Ho\A(QX,QY)$ a group homomorphism. The fact that the composition in $\Ho\A$ is bilinear easily follows from properties of left and right homotopies (see~\cite[\S1.2]{Hov99}).
\end{proof}

An easy corollary of the proposition demonstrates why we shall later consider the derived category $\Der(\A)$ as $\Ho\Cpx(\A)$ for a suitable model structure on $\Cpx(\A)$ rather than using the more classical algebraic approach via Ore localization of the homotopy category of complexes $\Htp(\A)$ (as in~\cite{Hap88,Verd96}). Expressing $\Ho\Cpx(\A)(X,Y)$ as an $\Ext^1$ group in $\Cpx(\A)$ will allow us to apply Proposition~\ref{prop:eklof} and argue using filtrations in $\Cpx(\A)$. This is something not readily available in the language of triangulated categories.

\begin{cor} \label{cor:hom_in_htpy}
Let $\A$ be a weakly idempotent complete exact category and let $(\Cof,\We,\Fib)$ encode an exact model structure on $\A$. Suppose that $X$ is cofibrant, $Y$ is fibrant and $0 \to \Omega Y \to P \to Y \to 0$ is a conflation with $\Omega Y \in \Fib$ and $P \in \Cof\cap\We$ (such a conflation always exists by the completeness of $(\Cof\cap\We, \Fib)$). Then naturally $\Ho\A(QX,QY) \cong \Ext^1_\A(X,\Omega Y)$.
\end{cor}

\begin{proof}
Since $\F$ is closed under extensions, we have $P \in \Cof\cap\We\cap\Fib$. Applying $\A(Z,-)$ with arbitrary $Z \in \A$ to $0 \to \Omega Y \to P \to Y \to 0$, we obtain an exact sequence
\[ \A(Z,P) \la \A(Z,Y) \la \Ext^1_\A(Z,\Omega Y) \la \Ext^1(Z,P) \]
If $Z \in \Cof\cap\We\cap\Fib$, then $\Ext^1_\A(Z,\Omega Y) = 0$. Therefore, a morphism $f\colon X \to Y$ factors through an object of $\Cof\cap\We\cap\Fib$ \iff it factors through the deflation $P \to Y$. If $Z = X$ is cofibrant, then $\Ext^1(X,P) = 0$ and we have
\[ \A(X,P) \la \A(X,Y) \la \Ext^1_\A(X,\Omega Y) \la 0. \]
Finally, $\Ho\A(QX,QY)$ is also isomorphic to the cokernel of $\A(X,P) \to \A(X,Y)$ by Proposition~\ref{prop:hom_in_htpy} and the above consideration, so that we have $\Ho\A(QX,QY) \cong \Ext^1_\A(X,\Omega Y)$.
\end{proof}

We conclude the section by discussing when the homotopy category of an exact model structure carries a natural triangulated structure (or in homotopy theoretic terminology, when the model structure is stable). A very mild sufficient condition is provided in~\cite[\S1]{Beck14}.

\begin{defn} \label{defn:hereditary}
A cotorsion pair $(\X,\Y)$ in an exact category $\A$ is called \emph{hereditary} if $\Ext^i_\A(X,Y) = 0$ for all $X \in\X$, $Y \in\Y$ and \emph{all} $i>0$. An exact model structure on $\A$ given by $(\Cof,\We,\Fib)$ is called a \emph{hereditary model structure} if the two complete cotorsion pairs $(\Cof,\We\cap\Fib)$ and $(\Cof\cap\We,\Fib)$ are hereditary.
\end{defn}

\begin{rem} \label{rem:hereditary}
The reader who is possibly confused by seemingly different definitions of hereditary cotorsion pairs in the literature is referred to~\cite[Corollary 1.1.12]{Beck14} and~\cite[Lemma 6.17]{St14} for rectification.
\end{rem}

\begin{prop} \label{prop:model-triang}
Suppose that $\A$ is a weakly idempotent complete exact category and $(\Cof,\We,\Fib)$ represents a hereditary model structure on $\A$. Then $\Ho\A$ carries a natural triangulated structure \st every conflation $0 \to X \stackrel{u}\to Y \stackrel{u}\to Z \to 0$ in $\A$ yields a triangle $QX \stackrel{Qu}\to QY \stackrel{Qu}\to QZ \stackrel{Qw}\to \Sigma(QX)$ in $\Ho\A$ and, moreover, every triangle arises from a conflation in this way.
\end{prop}

\begin{proof}
The proof of this fact for hereditary model structures on abelian categories from~\cite[Proposition 1.1.15]{Beck14} easily carries over to arbitrary hereditary model structures; see~\cite[Theorem 6.21]{St14}. The main point is that $\Cof\cap\Fib$ with the exact structure induced from $\A$ becomes a Frobenius exact category whose projective-injective objects are precisely those in $\Cof\cap\We\cap\Fib$. Then $\Ho\A$ can be canonically identified with the stable category of $\Cof\cap\Fib$ which is known to be triangulated by~\cite[Theorem I.2.6]{Hap88}.
\end{proof}

\section{Purity in finitely accessible categories}
\label{sec:purity}

\subsection{Pure exact structures}
\label{subsec:pure-ex}

Here we describe a class of exact structures which are prominent in this paper: pure exact structures on additive finitely accessible categories. Although finitely accessible categories may seem rather abstract, there is a handy representation result in terms of flat modules which we are going to use in our arguments.

\begin{defn} \label{defn:flat}
Let $\clS$ be skeletally small preadditive category. An $\clS$-module $F \in \Modu\clS$ is called \emph{flat} if it is a direct limit of finitely generated projective $\clS$-modules. We denote the full subcategory of $\Modu\clS$ consisting of all flat modules by $\Flat\clS$.
\end{defn}

For other characterizations and basic properties of flat modules we refer for example to~\cite[\S1.3]{CB94}. For the moment we only point out that $\Flat\clS$ is closed in $\Modu\clS$ under coproducts and direct limits. In particular, $\Flat\clS$ always has coproducts and direct limits as a category. Now we can explain the connection to finitely accessible categories.

\begin{prop}
\label{prop:rep-th-fa}
The assignments
$$ \clS \mapsto \Flat\clS \quad \textrm{ and } \quad \A \mapsto \fp\A $$
(recall Definition~\ref{defn:fin-acc}) give a bijective correspondence between
\begin{enumerate}
\item equivalence classes of skeletally small idempotent complete additive categories $\clS$, and
\item equivalence classes of additive finitely accessible categories $\A$.
\end{enumerate}
\end{prop}

\begin{proof}
See~\cite[\S 1.4]{CB94} or~\cite[Theorem 2.26]{AR94}.
\end{proof}

Since $\Flat\clS$ is well-known to be an extension closed subcategory of $\Modu\clS$, it follows from~\cite[Lemma 10.20]{Bu10} that an additive finitely accessible category always carries a natural exact structure. To be more precise, Proposition~\ref{prop:rep-th-fa} identifies $\A$ with $\Flat\clS$, and under this identification $\fp\A$ corresponds to $\proj\clS$, the full subcategory of finitely generated projective modules. Now $\Flat\clS$ carries the exact structure inherited from $\Modu\clS$: Conflations are simply the short exact sequences in $\Modu\clS$ whose all three terms are flat. In particular, conflations in $\Flat\clS$ are characterized by the property that $\Hom_\clS(P,-)$ sends them to short exact sequences of abelian groups for each $P \in \proj\clS$. Translating this characterization back to $\A$ leads to the following definition from~\cite{CB94}.
\begin{defn} \label{defn:pure-exact}
Let $\A$ be an additive finitely accessible category. A sequence $0 \to X \to Y \to Z \to 0$ in $\A$ is called \emph{pure exact} if the induced sequence of abelian groups
\[ 0 \la \A(F,X) \la \A(F,Y) \la \A(F,Z) \la 0 \]
is exact for each $F \in \fp\G$. The \emph{pure exact structure} on $\A$ is the one where conflations are precisely all pure exact sequences. Then we call the inflations \emph{pure monomorphisms} and the deflations \emph{pure epimorphisms}.
\end{defn}

A favorable fact is that both pure monomorphisms and epimorphisms admit a simple categorical characterization in $\A$.

\begin{lem} \label{lem:char-pure-mono-epi}
Let $\A$ be an additive finitely accessible category and $f\colon X \to Y$ be a morphism. Then
\begin{enumerate}
\item $f$ is a pure epimorphism \iff $\A(F,f)$ is surjective for each $F \in \fp\A$.
\item $f$ is a pure monomorphism \iff in each commutative square
\[
\begin{CD}
F @>{f'}>> G     \\
@V{u}VV @VV{v}V  \\
X @>{f}>>  Y
\end{CD}
\]
with $F,G \in \fp\A$, $u$ factors through $f'$.
\end{enumerate}
\end{lem}

\begin{proof}
In both cases we assume without loss of generality that $\A = \Flat\clS$ for $\clS$ as in Proposition~\ref{prop:rep-th-fa}.

(1) We only need to prove that $f$ has a kernel in $\A = \Flat\clS$. But as $\fp\A$ identifies with $\proj\clS$, the map $f$ is an epimorphism in $\Modu\clS$. Since $Y$ is flat, $f$ is known to be even a pure epimorphism in $\Modu\clS$ and so its kernel in $\Modu\clS$ is flat.

(2) If $f$ satisfies the factorization condition, then $f$ is necessarily a mono\-mor\-phism in $\Modu\clS$ by \cite[Proposition 2.29]{AR94}. Hence the factorization condition is equivalent to $f$ being a pure monomorphism in $\Modu\clS$ by \cite[Proposition 7.16]{JL89}. This is in turn equivalent to $f$ having a flat cokernel in $\Modu\clS$, which is none other than saying that $f$ is a pure monomorphism in $\A = \Flat\clS$.
\end{proof}

We also note that if $\A = \Flat\clS$ with the pure exact structure, then we have canonical natural isomorphisms $\Ext^i_\A(X,Y) \cong \Ext^i_\clS(X,Y)$ for each $X,Y\in\A$ and $i \ge 0$ by (the proof of)~\cite[Lemma 14]{OpSt12}.

As syzygies of flat modules are known to be flat, $\Flat\clS$ has always enough projective objects as an exact category (see~\cite[Definition~11.1]{Bu10}). The projectives of $\Flat\clS$ are of course the usual projective $\clS$-modules. Translating this back to the setting of abstract additive finitely accessible categories, we see that $\A$ with the pure exact structure has enough projectives, which we call the \emph{pure projective} objects of $\A$. The class of all pure projective objects will be denoted by $\PProj\A$ and consists precisely of summands of arbitrary direct sums of finitely presentable objects.

A fact which may be less well known is that the pure exact structure also admits enough injectives. If $\A$ has products, this has been proved in~\cite[\S3.3]{CB94}. In general this is a consequence of the solution to the Flat Cover Conjecture~\cite{BEE01}.

\begin{defn} \label{defn:cot-obj}
Let $\clS$ be skeletally small preadditive category. We will say that an $\clS$-module $Y$ is a \emph{cotorsion} module if $\Ext^1_\clS(F,Y) = 0$ for all flat $\clS$-modules $F$. The class of all cotorsion modules will be denoted by $\Cot\clS$. 
\end{defn}

\begin{rem} \label{rem:cot-obj}
The terminology comes from the fact that $\Cot\Z$ was a part of the first cotorsion pair ever considered (see~\cite{Sal79}).
\end{rem}

\begin{prop} \label{prop:FCC}
Let $\clS$ be a skeletally small preadditive category. Then $(\Flat\clS,\Cot\clS)$ is a functorially complete hereditary cotorsion pair in $\Modu\clS$. In particular, if we denote $\A = \Flat\clS$ and $\I = \Flat\clS\cap\Cot\clS$, then $(\A,\I)$ is a functorially complete cotorsion pair in $\A$ and $\A$ has enough injectives as an exact category.
\end{prop}

\begin{proof}
As mentioned, the key ingredient is the solution to the Flat Cover Conjecture from~\cite{BEE01}.
The pair $(\Flat\clS,\Cot\clS)$ is a complete cotorsion pair in $\Modu\clS$ since the proof in \cite[\S2]{BEE01} for modules over rings directly generalizes to modules over small categories.
The functoriality of approximation sequences follows from Quillen's small object argument (see~\cite[Theorem 5.16]{St14}) and the heredity of $(\Flat\clS,\Cot\clS)$ follows from the fact that flat modules are closed under syzygies and standard dimension shifting.
For the second claim, note that $\Ext^1_\A(A,I) \cong \Ext^1_\clS(A,I) = 0$ for each $A \in \A$ and $I \in \I$. Moreover, every $A \in \A$ admits functorial approximation sequences \wrt $(\A,\I)$ and both $\A,\I$ are closed under retracts. It easily follows that $(\A,\I)$ is a functorially complete cotorsion pair in $\A$. Since the objects of $\I$ are by definition precisely the injective objects in $\A$, we have shown that $\A$ has enough injectives.
\end{proof}

\begin{defn} \label{defn:pinj}
If $\A$ is a finitely accessible additive category, an object injective with respect to the pure exact structure on $\A$ will be called \emph{pure injective}. We denote the class of all pure injective objects by $\PInj\A$.
\end{defn}

\subsection{Categories of complexes}
\label{subsec:complexes}

Later we will also need to work with the category of complexes over an additive finitely accessible category $\A$ as well as with various related categories. The reason why we discuss this here is that the category of complexes carries at least two or three natural exact structures, and we need to distinguish between these carefully to prevent confusion.

We will use the cohomological notation for complexes:
\[ X\colon \quad \cdots \la X^{-1} \stackrel{d^{-1}_X}\la X^0 \stackrel{d_X^0}\la X^1 \stackrel{d_X^1}\la X^2 \la \cdots. \]
The category of all complexes with usual cochain complex homomorphisms will be denoted by~$\Cpx(\A)$, and the additive quotient by the null-homotopic maps will be denoted by $\Htp(\A)$. The suspension functor, which shifts each $X \in \Cpx(\A)$ by one degree to the left and changes the signs of all differentials, will be denoted by $\Sigma\colon \Cpx(\A) \to \Cpx(\A)$. We will view $\A$ as a full subcategory of $\Cpx(\A)$ via the embedding which sends $A \in \A$ to the complex
\[ \cdots \la 0 \la 0 \la A \la 0 \la 0 \la \cdots \]
with $A$ in degree zero.

We will be interested in two distinct exact structures on $\Cpx(\A)$ and we will give them names for future reference.

\begin{defn} \label{defn:exact-cpx}
The exact structure on $\Cpx(\A)$ whose conflations are defined by the property that they are pure exact in $\A$ in each degree, will be called the \emph{induced pure exact} structure on $\Cpx(\A)$. By somewhat abusing the term, we will call such conflations \emph{pure exact sequences} of complexes. The inflations and deflations will be called \emph{pure monomorphisms} and \emph{pure epimorphisms} of complexes, respectively. We denote the Ext groups \wrt this exact structure by $\Ext^i_{\Cpx(\A)}$.

The second exact structure which is often considered on $\Cpx(\A)$ is the \emph{degree wise split} exact structure, whose conflations are those sequences of complexes which are split exact in each degree. We denote the Ext groups \wrt this exact structure by $\Ext^i_{\Cpx(\A),\dw}$.
\end{defn}

\begin{rem} \label{rem:pe-induced}
The category $\Cpx(\A)$ is always a $\kappa$-accessible category for some $\kappa$ (Definition~\ref{defn:fin-acc}) by \cite[Proposition 4.5]{SaSt11}. Although this does not guarantee that $\Cpx(\A)$ is finitely accessible, this is often the case (e.g.\ if $\A$ is abelian, see Remark~\ref{rem:pure-acyclic-Groth} below). Thus, $\Cpx(\A)$ may carry its own inherent pure exact structure, which is usually different from the induced pure exact structure and which we never shall consider here. The word \emph{induced} is to stress this difference.
\end{rem}

We will also consider the usual $\Cpx(\Ab)$-enriched version of $\Cpx(\A)$. For $X,Y \in \Cpx(\A)$ we denote by $\HOM_\A(X,Y)$ the complex of abelian groups \st
\[ \HOM_\A(X,Y)^n = \prod_{i \in \Z} \A(X^i,Y^{i+n}) \]
and the differential for $f \in \HOM_\A(X,Y)^n$ is given by the graded commutator $d(f) = d_Y \circ f - (-1)^n f \circ d_X$. Then clearly $H^n \HOM_\A(X,Y) \cong \Htp(\A)(X,\Sigma^n Y)$. The complexes over $\A$ with such homomorphism complexes form a dg-category~\cite{Kel94}, but we will not need this fact. For our purpose it suffices to note the following.

\begin{samepage}
\begin{lem} \label{lem:HOM} ~
\begin{enumerate}
\item $\HOM_\A(-,Y)\colon \Cpx(\A)\op \to \Cpx(\Ab)$ sends pure exact sequences of complexes to left exact sequences in $\Cpx(\Ab)$ (i.e.\ $\HOM_\A(-,Y)$ is left exact).
\item If all components of $Y$ are pure injective, then $\HOM_\A(-,Y)$ sends pure exact sequences to exact sequences (i.e.\ $\HOM_\A(-,Y)$ is exact in this case).
\end{enumerate}
\end{lem}
\end{samepage}

\begin{proof}
This is straightforward from the left exactness of $\A(-,Y^i)$ and its exactness if $Y^i$ is pure injective.
\end{proof}

Another standard category related to $\Cpx(\A)$ is the derived category $\Der(\A)$ of $\A$ (see \cite[\S10]{Bu10}). For the sake of completeness we state the relevant definitions here explicitly. We shall discuss this category more in detail in Section~\ref{sec:pure-derived}.

\begin{defn}[\cite{ChHo02,Kr12}] \label{defn:pure-acyclic}
A complex $X \in \Cpx(\A)$ is \emph{pure acyclic} if each differential $d^{i-1}\colon X^{i-1} \to X^{i}$ factors as $X^{i-1} \to Z^{i}(X) \to X^i$ in such a way that each of the induced sequences $0 \to Z^i(X) \to X^i \to Z^{i+1}(X) \to 0$ is pure exact in $\A$. A map of complexes $f\colon X \to Y$ is a \emph{pure quasi-isomorphism} if the mapping cone of $f$ is pure acyclic. The \emph{pure derived category} $\Der(\A)$ of $\A$ is defined as the localization of $\Cpx(\A)[\we\inv]$ where $\we$ is the class of all pure quasi-isomorphisms.
\end{defn}

\begin{rem} \label{rem:pure-acyclic-Groth}
Suppose that $\G$ is a finitely accessible abelian category. Then it is automatically locally finitely presentable and it is also a Grothendieck category, see~\cite[\S2.4]{CB94}. Moreover, $\Cpx(\G)$ again a locally finitely presentable Grothendieck category by~\cite[Corollary 1.54]{AR94} (see also~\cite[Lemma 1.1]{St13}). Since $\fp\G$ is a generating set, one easily checks in this case that $X$ is pure acyclic \iff $\Htp(\G)(\Sigma^i A,X) = 0$ for each $A \in \fp\G$ and $i \in \Z$.
\end{rem}

We conclude the section by two useful lemmas. The first one provides a convenient ``generating set'' for $\Cpac(\A)$ in terms of direct limits.

\begin{lem} \label{lem:gens-for-CA}
Let $\A$ be an additive finitely accessible category. A complex $X$ belongs to $\Cpac(\A)$ \iff it is a direct limit of bounded contractible complexes over $\fp\A$.
\end{lem}

\begin{proof}
By Proposition~\ref{prop:rep-th-fa}, we can assume that $\A = \Flat\clS$ and $\fp\A$ is identified with the full subcategory of finitely generated projective $\clS$-modules. Then the statement was proved in~\cite[Theorem 2.4]{EGR98}. For an alternative proof we refer to~\cite[Facts 2.4 and Theorem 8.6]{Nee08}.
\end{proof}

The second one gives a certain adjunction for the $\Ext$-functor between $\A$ and $\Cpx(\A)$. Given $X \in \A$ and $n \in \Z$, we denote by $D^n(X)$ the complex
\begin{equation} \label{eq:discs}
D^n(X)\colon \quad \cdots \la 0 \la 0 \la X \stackrel{1_X}\la X \la 0 \la 0 \la \cdots 
\end{equation}
with $X$ in cohomological degrees $n$ and $n+1$.

\begin{lem} \label{lem:ext-adj}
Let $\A$ be an additive finitely accessible category. Given $X \in \A$, $n \in \Z$ and $Y \in \Cpx(\A)$, there are natural isomorphisms
\begin{enumerate}
\item $\Ext^1_\A(X,Y^n) \cong \Ext^1_{\Cpx(\A)}(D^n(X),Y)$ and
\item $\Ext^1_\A(Y^n,X) \cong \Ext^1_{\Cpx(\A)}(Y,D^{n-1}(X))$.
\end{enumerate}
\end{lem}

\begin{proof}
Exactly the same proof as for~\cite[Lemma 3.1(5) and (6)]{Gil04} applies to arbitrary exact categories.
\end{proof}

\section{Models for the pure derived category}
\label{sec:pure-derived}

This section is devoted to studying the pure derived category $\Der(\A)$ of an additive finitely accessible category $\A$. The results will be useful in Sections~\ref{sec:coderived} and~\ref{sec:singularity} for a conceptual explanation of the role of purity \wrt coderived and singularity categories.

As a matter of the fact, an important step in this context which was taken by Neeman~\cite{Nee08}, was motivated by coderived categories. What he proved can be stated in our terminology as follows (see Theorem~\ref{thm:cotorsion-pac} below): If $\A$ is an additive finitely accessible category, then $\big(\Cpx(\PProj\A), \Cpac(\A)\big)$ is a complete cotorsion pair in $\Cpx(\A)$ \wrt the induced pure exact structure. In particular, a pure acyclic complex in $\Cpx(\A)$ with pure projective components is always contractible. Here we shall prove dual results: $\big(\Cpx(\A), \Cpac(\PInj\A)\big)$ is also a complete cotorsion pair in $\Cpx(\A)$, and a pure acyclic complex with pure injective components also must be contractible. This has implications for model structures for $\Der(\A)$ which, in more classical algebraic terminology, can be restated in terms of a recollement.

We start the discussion with an easy but crucial observation.

\begin{lem} \label{lem:inj-2-of-3}
Let $\A$ be an additive finitely accessible category and $\Y \subseteq \Cpx(\PInj\A)$ be a class complexes \st $\Y = \Sigma\Y$. Then the class
\[ \X = \{ X \in \Cpx(\A) \mid \Ext^1_{\Cpx(\A)}(X,Y) = 0 \textrm{ for all } Y \in \Y \} \]
has the 2-out-of-3 property for pure exact sequences. That is, if $0 \to A \to B \to C \to 0$ is a conflation in $\Cpx(\A)$ with the induced pure exact structure (Definition~\ref{defn:exact-cpx}) and two of $A,B,C$ belong to $\X$, so does the third.
\end{lem}

\begin{proof}
Consider $Y \in \Cpx(\PInj\A)$. Applying the functor $\HOM_\A(-,Y)$ to $0 \to A \to B \to C \to 0$, we get a short exact sequence
\[ 0 \la \HOM_\A(C,Y) \la \HOM_\A(B,Y) \la \HOM_\A(A,Y) \la 0 \]
of complexes of abelian groups by Lemma~\ref{lem:HOM}. Recall that for each $X \in \Cpx(\A)$ and $Y \in \Cpx(\PInj\A)$ we have
\begin{multline*}
H^n(\HOM_\A(X,Y)) \cong \Htp(\A)(X,\Sigma^n Y) \cong \\
\cong \Ext^1_{\Cpx(\A),\dw}(X,\Sigma^{n-1} Y) \cong \Ext^1_{\Cpx(\A)}(X,\Sigma^{n-1} Y),
\end{multline*}
for all $n \in \Z$. The first two isomorphisms are well known, while the last one follows from the fact that $Y$ has pure injective components, so every degree wise pure exact sequence is degree wise split. Thus, we obtain a long exact sequence of groups
\begin{multline*}
\cdots \la \Ext^1_{\Cpx(\A)}(A,\Sigma^{n-1} Y) \la \Ext^1_{\Cpx(\A)}(C,\Sigma^{n} Y) \la \\
\la \Ext^1_{\Cpx(\A)}(B,\Sigma^{n} Y) \la \Ext^1_{\Cpx(\A)}(A,\Sigma^{n} Y) \la \\
\la \Ext^1_{\Cpx(\A)}(C,\Sigma^{n+1} Y) \la \Ext^1_{\Cpx(\A)}(B,\Sigma^{n+1} Y) \la \cdots
\end{multline*}
The 2-out-of-3 property follows immediately.
\end{proof}

It follows directly from Definition~\ref{defn:pure-exact} that a direct limit of pure exact sequences in $\A$ is again pure exact. Thus, it makes sense to speak of \emph{direct unions} in $\A$; these are the colimits of direct systems $(X_i, f_{ij} \mid i<j;\, i,j \in I)$ with all $f_{ij}$ pure monomorphisms. Similarly, since $\Cpx(\A)$ has arbitrary direct limits which are computed component wise, we can speak of direct unions in $\Cpx(\A)$ \wrt the induced pure exact structure. The following observation then essentially follows from the previous lemma and Proposition~\ref{prop:eklof}.

\begin{lem} \label{lem:inj-union}
Let $\A$ be additive finitely accessible and $\Y \subseteq \Cpx(\PInj\A)$ be \st $\Y = \Sigma\Y$. Then
$\X = \{ X \in \Cpx(\A) \mid \Ext^1_{\Cpx(\A)}(X,Y) = 0 \textrm{ for all } Y \in \Y \}$
is closed under direct unions in $\Cpx(\A)$.
\end{lem}

\begin{proof}
Let $(X_i, f_{ij} \mid i<j;\, i,j \in I)$ be a direct system in $\X$ \st all $f_{ij}$ are pure monomorphisms in the sense of Definition~\ref{defn:exact-cpx}. Using a standard reduction as in~\cite[Corollary 1.7]{AR94}, we can without loss of generality assume that the direct system is a $\lambda$-sequence for an infinite cardinal number $\lambda$, so in particular of the form $(X_\alpha, f_{\alpha\beta} \mid \alpha<\beta<\lambda)$.

Note that then $\Coker f_{\alpha,\alpha+1} \in \X$ by Lemma~\ref{lem:inj-2-of-3}, so that we almost have an $\X$-filtration in the sense of Definition~\ref{defn:filtr}. In fact we can adjust the direct system to obtain an $\X$-filtration $(X_\alpha, f_{\alpha\beta} \mid \alpha<\beta\le\lambda)$: we add the zero object at the beginning and construct $X = X_\lambda$ as $\varinjlim_{\alpha<\lambda} X_\alpha$. Now $X \in \X$ by Proposition~\ref{prop:eklof}, but $X$ is none other than the direct union of $(X_\alpha, f_{\alpha\beta} \mid \alpha<\beta<\lambda)$.
\end{proof}

As a consequence of the previous lemmas, we obtain a key result.

\begin{prop} \label{prop:inj-pac}
Let $\A$ be an additive finitely accessible category. If $\Y \subseteq \Cpx(\PInj\A)$ is a class complexes \st $\Y = \Sigma\Y$, then the class
\[ \X = \{ X \in \Cpx(\A) \mid \Ext^1_{\Cpx(\A)}(X,Y) = 0 \textrm{ for all } Y \in \Y \} \]
is closed under direct limits in $\Cpx(\A)$.
\end{prop}

\begin{proof}
Suppose that $(X_i, f_{ij} \mid i<j;\, i,j \in I)$ is a direct system in $\X$ with $X = \varinjlim X_i$ in $\Cpx(\A)$, and consider the morphism $p\colon \coprod_i X_i \to X$ with the colimit maps as components. First observe that $p$ has a kernel in $\Cpx(\A)$. Indeed, if we identify $\A = \Flat\clS$ according to Proposition~\ref{prop:rep-th-fa}, then $p$ is a surjective map in $\Cpx(\Modu\clS)$. Therefore, the kernel of $p$ belongs to $\Cpx(\Flat\clS)$ as the class of flat modules is well known to be closed under kernels of epimorphisms in $\Modu\clS$ (see e.g.~\cite[Proposition 2.2]{Len83}).

Thus, we have a pure exact sequence of the following form
\[ 0 \la K \stackrel{\ell}\la \coprod_i X_i \la X \la 0 \]
in $\Cpx(\A)$. Instead of proving directly that $X = \varinjlim X_i \in \X$, thanks to Lemma~\ref{lem:inj-2-of-3} it suffices to prove that $K \in \X$. This is what we aim to do next by expressing $K$ as a direct union of objects from $\X$.

Almost by definition, $p\colon \coprod_i X_i \to X$ is the cokernel of a morphism $h\colon \coprod_{i<j} X_{ij} \to \coprod_i X_i$ given as follows. We put $X_{ij} = X_i$ for each $i<j$ and the component $h_{ij}\colon X_{ij} \to \coprod_i X_i$ of $h$ is given as the composition of $(^{\;1_{X_i}}_{-f_{ij}})\colon X_{ij} \to X_i \oplus X_j$ with the coproduct inclusion $X_i \oplus X_j \to \coprod_i X_i$.

Suppose now that $J \subseteq I$ is a finite directed subposet, or in other words, $J$ is a subset of $I$ with a unique maximal element with respect to the ordering on $I$. Identifying $\Cpx(\A)$ with $\Cpx(\Flat\clS)$ again, we shall denote by $K_J$ the image in $\Cpx(\Modu\clS)$ of the map $(h_{ij})\colon \coprod_{i} X_{ij} \to K$, where $j = \max J$ and $i$ runs through all non-maximal elements of $J$. A direct computation shows that the composition
\[ \coprod_{i \in J \setminus\{\max J\}} X_{ij} \stackrel{(h_{ij})}\la K \stackrel{\ell} \la \coprod_{i \in I} X_i \la \coprod_{i \in J \setminus\{\max J\}} X_i, \]
where the last arrow is the canonical projection, is the identity. In particular, $K_J$ is a summand of $K$ and as such belongs to $\Cpx(\A)$, and in fact even $K_J \cong \coprod_{i \in J \setminus\{\max J\}} X_i \in \X$.
Another simple observation is that if $J \subseteq J'$ are two finite directed subposets in $I$, then $K_J \subseteq K_{J'}$ (as complexes of $\clS$-modules). That is, $K$ is the direct union of $(K_J \mid J \textrm{ finite directed subposet of } I)$ both in $\Cpx(\Modu\clS)$ and in $\Cpx(\A)$. As all $K_J \in \X$ for every $J$, we have $K \in \X$ by Lemma~\ref{lem:inj-union}. Thus, as mentioned above, $X = \varinjlim_i X_i \in \X$ by Lemma~\ref{lem:inj-2-of-3}.
\end{proof}

Now we are in a position to prove the existence of the aforementioned complete cotorsion pairs.

\begin{thm} \label{thm:cotorsion-pac}
Let $\A$ be an additive finitely accessible category. Then
\[ \big(\Cpx(\PProj\A), \Cpac(\A)\big) \qquad \textrm{and} \qquad \big(\Cpac(\A), \Cpx(\PInj\A)\big) \]
are functorially complete hereditary cotorsion pairs in $\Cpx(\A)$ with the induced pure exact structure.
\end{thm}

\begin{proof}
As before, we can identify $\A$ with $\Flat\clS$ for $\clS = \fp\A$. Then $\Cpx(\PProj\A) = \Cpx(\Proj\clS)$ is deconstructible in $\Cpx(\Modu\clS)$ (with the abelian exact structure) by~\cite[Theorem 4.2(1)]{St13}. It follows that $\Cpx(\PProj\A)$ is also deconstructible in $\Cpx(\A) = \Cpx(\Flat\clS)$ with the induced pure exact structure. In other words, there is a small set $\U \subseteq \Cpx(\PProj\A)$ \st $\Cpx(\PProj\A) = \Filt(\U)$ in $\Cpx(\A)$. Using the same notation as in~\eqref{eq:discs} and Lemma~\ref{lem:ext-adj}, we can without loss of generality assume that up to isomorphism all complexes of the form
\[
D^n(F)\colon \quad \cdots \la 0 \la 0 \la F \stackrel{1_F}\la F \la 0 \la 0 \la \cdots,
\]
with $F \in \fp\A$ and $n \in \Z$, belong to $\U$. Since these complexes form a generating set for $\Cpx(\A)$, we can apply Proposition~\ref{prop:complete-cotor-A} and this proves the existence of a functorially complete cotorsion pair $\big(\Cpx(\PProj\A), \Y\big)$ in $\Cpx(\A)$. Since $\Cpx(\PProj\A)$ is closed under taking syzygies in $\Cpx(\A)$, the cotorsion pair is hereditary by dimension shifting.

We claim that $\Y = \Cpac(\A)$. To this end, we learn from~\cite[Theorem 8.6]{Nee08} that $Y \in \Cpx(\A)$ belongs to $\Cpac(\A)$ \iff $\Htp(\A)(X,Y) = 0$ for all $X \in \Cpx(\PProj\A)$. But as in the proof of Lemma~\ref{lem:inj-2-of-3} we have isomorphisms
\[ \Htp(\A)(X,Y) \cong \Ext^1_{\Cpx(\A),\dw}(\Sigma X,Y) \cong \Ext^1_{\Cpx(\A)}(\Sigma X,Y), \]
proving the claim.

Regarding the other cotorsion pair, $\Cpac(\A) = \Cpac(\Flat\clS)$ is deconstructible in $\Cpx(\Modu\clS)$ by~\cite[Theorem 4.2(2)]{St13}. Hence $\Cpac(\A)$ is also deconstructible in $\Cpx(\A)$. Again, $\Cpac(\A)$ contains all the complexes $D^n(F)$ with $F \in \fp\A$ and $n \in \Z$, and is closed under taking syzygies in $\Cpx(\A)$. Hence we have a functorially complete hereditary cotorsion pair $\big(\Cpac(\A), \Y'\big)$ in $\Cpx(\A)$. 

It remains to show that $\Y' = \Cpx(\PInj\A)$. As $\Cpac(\A)$ in fact contains all $D^n(X)$ for arbitrary $X \in \A$ and $n \in \Z$, it follows from Lemma~\ref{lem:ext-adj}(1) that $\Y' \subseteq \Cpx(\PInj\A)$. On the other hand, the class
\[ \X' = \{ X \in \Cpx(\A) \mid \Ext^1_{\Cpx(\A)}(X,Y) = 0 \textrm{ for all } Y \in \Cpx(\PInj\A) \} \]
contains all bounded contractible complexes with finitely presented components and is closed under direct limits by Proposition~\ref{prop:inj-pac}. Therefore, $\X' \supseteq \Cpac(\A)$ by Lemma~\ref{lem:gens-for-CA}. Hence $\Y' \supseteq \Cpx(\PInj\A)$, concluding the proof.
\end{proof}

In the rest of the section we state several corollaries which illustrate the power of the theorem. First of all, we obtain the following rather non-obvious criterion for contractibility of a complex with pure-injective components.

\begin{cor} \label{cor:pure-acyclic-pinj}
Suppose that $X \in \Cpx(\A)$ is a pure acyclic complex. If all components of $X$ are pure projective or if all components of $X$ are pure injective, then $X$ is contractible.
\end{cor}

\begin{proof}
We will only prove the case with $X \in \Cpx(\PInj\A)$. The other case is similar and it has been obtained before by using a different argumentation; see~\cite[Remark 2.15]{Nee08}. So let $X \in \Cpx(\PInj\A) \cap \Cpac(\A)$. Then $\Ext^1(\Sigma X,X) = 0$ by Theorem~\ref{thm:cotorsion-pac}. Again we have standard isomorphisms
\[ \Ext^1_{\Cpx(\A)}(\Sigma X,X) \cong \Ext^1_{\Cpx(\A),\dw}(\Sigma X,X) \cong \Htp(\A)(X,X). \]
Thus, $1_X$ is null-homotopic, or equivalently $X$ is contractible.
\end{proof}

Second, it follows that each pure acyclic complex over $\A$ can be transfinitely built from complexes as in \eqref{eq:discs} above Lemma~\ref{lem:ext-adj}. If we replaced $\A$ by a category of quasi-coherent sheaves with the abelian exact structure, this would be interpreted as a regularity condition. Hence additive finitely accessible categories behave as regular in this sense, regardless of whether their global dimension is finite or not.

\begin{cor} \label{cor:pure-acyclic-cellular}
There is a small set $\U \subseteq \A$ \st each pure acyclic complex $X \in \Cpac(\A)$ is a retract of a $\{D^n(U)\mid U \in \U \textrm{ and } n \in \Z\}$-filtered complex.
\end{cor}

\begin{proof}
Assume again that $\A = \Flat\clS$. It is well known that the cotorsion pair $(\Flat\clS,\Cot\clS)$ in $\Modu\clS$ is generated by a small set $\U$ in the sense of Proposition~\ref{prop:complete-cotor-A}; see~\cite[Proposition 2]{BEE01}. Of course, we can without loss of generality assume that all isotypes of finitely presentable objects are represented in $\U$. It follows that
\[ \PInj\A = \{ Y \in \A \mid \Ext^1_\A(U,Y) = 0 \textrm{ for each } U \in \U \} \]
and, in view of Lemma~\ref{lem:ext-adj}, also that
\[ \Cpx(\PInj\A) = \{ Y' \in \Cpx(\A) \mid \Ext^1_{\Cpx(\A)}(D^n(U),Y') = 0 \textrm{ for each } U \in \U, n \in \Z \}. \]
Hence the conclusion is a consequence of Proposition~\ref{prop:complete-cotor-A} applied to $(\X,\Y) = \big(\Cpac(\A), \Cpx(\PInj\A)\big)$.
\end{proof}

In view of Proposition~\ref{prop:exact-model}, we also have an explicit description of two exact and hereditary model structures on $\Cpx(\A)$ whose homotopy category is the pure derived category $\Der(\A)$.
Following~\cite[1.1.10]{Beck14}, we call them the \emph{projective} and the \emph{injective model structure} for $\Der(\A)$, respectively.
The projective model structure has been described also in~\cite[Theorem 7.17]{St14}.

\begin{cor} \label{cor:models}
The triples
\[ (\Cof,\We,\Fib) = \big(\Cpx(\PProj\A),\Cpac(\A),\Cpx(\A)\big) \]
and
\[ (\Cof,\We,\Fib) = \big(\Cpx(\A),\Cpac(\A),\Cpx(\PInj\A)\big). \]
of classes of $\Cpx(\A)$ represent hereditary model structures on $\Cpx(\A)$. In both cases, the homotopy category is $\Der(\A)$.
\end{cor}

\begin{proof}
By Corollary~\ref{cor:pure-acyclic-pinj}, $\Cpx(\PProj\A) \cap \Cpac(\A)$ is the class of contractible complexes with pure projective components, or in other words the class of projective objects of $\Cpx(A)$ with the induced pure exact structure. Hence in the first case $(\Cof\cap\We,\Fib)$ is a (trivial) functorially complete cotorsion pair, and $(\Cof,\We,\Fib)$ satisfies the assumptions of Proposition~\ref{prop:exact-model}. The second case is similar. In this situation, $\Cpac(\A) \cap \Cpx(\PInj\A)$ is the class of injective objects of $\Cpx(\A)$ and, in view of the proof of Corollary~\ref{cor:pure-acyclic-cellular}, $(\Cof,\We\cap\Fib)$ is functorially complete in $\Cpx(\A)$ by Proposition~\ref{prop:complete-cotor-A}. In both cases, the trivial objects are precisely the pure acyclic complexes. It is standard to check using the description from Proposition~\ref{prop:exact-model} that weak equivalences are none other than pure quasi-isomorphisms, and so $\Ho\Cpx(\A) = \Der(\A)$.
\end{proof}

In more classical algebraic terms, the latter corollary can be interpreted as the existence of a certain recollement. In fact, the existence of the recollement has been established in Neeman's \cite[Theorem 3.1]{Nee10}, based on results in~\cite{BEIJR12}. The added value here is that we explicitly compute the essential image of $j_*$.

\begin{cor} \label{lem:recollement-pure}
Denote by $\Hpac(\A)$ the homotopy category of pure acyclic complexes over $\A$ and by $i_*\colon \Hpac(\A) \to \Htp(\A)$ the inclusion. Then there is a recollement
\[
\xymatrix@1{
\Hpac(\A) \ar[rr]|(.45)\hole|{i_*} &&
\Htp(\A) \ar@/_3ex/[ll]_{i^*} \ar@/^3ex/[ll]^{i^!} \ar[rr]|(.45)\hole|{j^*} &&
\Der(\A) \ar@/_3ex/[ll]_{j_!} \ar@/^3ex/[ll]^{j_*}
}.
\]
Moreover, the essential image of $j_!$ is $\Htp(\PProj\A)$ and the essential image of $j_*$ is $\Htp(\PInj\A)$.
\end{cor}

\begin{proof}
Applying \cite[Corollary 3.12]{SaSt11} to $\Cpx(\A)$ with the degree wise split and induced pure exact structures, we infer that the pairs of classes 
\[ \big(\Htp(\PProj\A), \Hpac(\A)\big) \qquad \textrm{and} \qquad \big(\Hpac(\A), \Htp(\PInj\A)\big) \]
are so-called Bousfield localizing pairs in $\Htp(\A)$. That is, the classes in each pair are mutually orthogonal \wrt the Hom-functor, they are closed under suspensions and desuspensions, and for each $A \in \A$ there are triangles, unique up to a unique isomorphism, of the form
\[ X^A \la A \la Y_A \la \Sigma X^A \qquad \textrm{and} \qquad Y^A \la A \la Z_A \la \Sigma Y^A  \]
with $X^A \in \Htp(\PProj\A)$, $Y_A,Y^A \in \Hpac(\A)$ and $Z_A \in \Htp(\PInj\A)$.
The conclusion is now a standard application of abstract results concerning Bousfield localization; see~\cite[Ch. 9]{Nee01} or \cite[\S4]{Kr10}.
\end{proof}

Finally, there is a puzzling consequence of Theorem~\ref{thm:cotorsion-pac} for module categories. As mentioned in the introduction, modules as in the coming corollary for instance inevitably appear in Gorenstein homological algebra. An application of the corollary will be given later in Proposition~\ref{prop:ding-chen}.

\begin{cor} \label{cor:gor-inj}
Let $\clS$ be a small additive category and let
\[ \cdots \la I_2 \la I_1 \la I_0 \la X \la 0 \]
be an exact sequence of $\clS$-modules with all $I_n$, $n \ge 0$, injective modules. Then $X$ is cotorsion (recall Definition~\ref{defn:cot-obj}). In fact, even $\Ext^1_\clS(F,X) = 0$ whenever $F$ has finite flat dimension.
\end{cor}

\begin{proof}
Let us extend the exact sequence to an acyclic complex of $\clS$-modules
\[ I\colon \quad \la I_2 \la I_1 \la I_0 \la I^0 \la I^1 \la \cdots \]
with injective components and $I^0$ in degree $0$. Our task is to prove that $\Htp(\Modu\clS)(F,\Sigma I) = 0$ for any fixed flat $\clS$-module $F$. To this end, consider the class $\Cac(\InjR)$ of all acyclic complexes with injective components. Then
\[ \X = \{ X \in \Cpx(\Modu\clS) \mid \Ext^1_{\Cpx(\Modu\clS)}(X,Y) = 0 \textrm{ for all } Y \in \Cac(\InjR) \} \]
is closed under direct limits by Proposition~\ref{prop:inj-pac}. Note that since every extension of $Y$ by $X$ is necessarily degree wise split, it does not matter whether we consider the abelian or the induced pure exact structure on $\Cpx(\Modu\clS)$. Since clearly $\Ext^1_{\Cpx(\Modu\clS)}(P,Y) \cong \Htp(\Modu\clS)(P,\Sigma Y) = 0$ for every finitely generated projective module $P$ and for each $Y \in \Cac(\InjR)$, it follows that also $F \in \X$. Hence $\Htp(\Modu\clS)(F,\Sigma I) \cong \Ext^1_{\Cpx(\Modu\clS)}(F,I) = 0$, as required.

The last part follows by easy dimension shifting, since for each $n\ge 0$ also $Z_n = \Coker(I_{n+1} \to I_n)$ is cotorsion and $\Ext^1_\clS(F,X) \cong  \Ext^n_\clS(F,Z_n) \cong \Ext^1_\clS(\Omega^{n-1}(F),Z_n)$.
\end{proof}

\begin{rem} \label{rem:gor-inj}
If $F \in \Modu\clS$ is an $\aleph_n$-presentable flat module, then the projective dimension of $F$ is at most $n$ by \cite[Theorem 2.45]{Os73} or \cite[Corollary 3.13]{Sim77} and hence $\Ext^1_\clS(F,X) = 0$ by dimension shifting. For a general flat module $F$ we are not aware of a more direct argument than the one presented above.
\end{rem}

\begin{rem} \label{rem:vNreg}
Suppose that $R$ is a von Neumann regular ring. Then every $R$-module is flat and cotorsion modules are precisely the injective ones. Corollary~\ref{cor:gor-inj} says that there cannot be any non-trivial cosyzygies of infinite order. If $R$ is of infinite global dimension, there exist, however, non-trivial cosyzygies of arbitrarily large finite order.
\end{rem}

\section{Coderived categories}
\label{sec:coderived}

Given a locally finitely presentable Grothendieck category $\G$, we can consider the thick subcategory $\C$ of $\Der(\G)$, the usual derived category of $\G$, generated by $\fp\G$. If $\G$ is \emph{locally coherent}, which by definition means that $\fp\G$ is abelian \cite[\S2.4]{CB94}, then $\C$ coincides with $\Db(\fp\G)$, the bounded derived category of $\fp\G$; see~\cite[Theorem 15.3.1]{KaSch06}.
It is well known that the objects of $\C$ may not be compact in $\Der(\G)$. As noticed by Krause in~\cite{Kr05}, for $\G$ locally noetherian one needs to enlarge $\Der(\G)$ to $\Htp(\Inj\G)$ and there $\C$ identifies with the compact objects. As $\Htp(\Inj\G)$ coincides by~\cite[Theorem 3.7]{Pos11} with what Positselski calls the \emph{coderived} category of $\G$, we will stick to this terminology here.

Our aim is to show that the local noetherianness can be replaced by local coherence. The proof is, however, much more involved for two reasons. First, we will rely on the non-trivial Corollary~\ref{cor:pure-acyclic-pinj}. Second, it is difficult to deal with coproducts in $\Htp(\Inj\G)$ directly, as is very well illustrated in~\cite{Nee14}. Hence, we will need an fp-injective model for the coderived category, allowing us to prove a version of~\cite[Conjecture 1]{Kr12} along the way. We will also illustrate in Example~\ref{expl:non-coh-fail} that the assumption of local coherence cannot be  easily dropped, not even for module categories.

\subsection{The derived category of fp-injectives}
\label{subsec:fpinj}

Our first aim is to better understand the derived category of the exact category of fp-injective objects.

\begin{defn} \label{defn:fp-inj}
Let $\G$ be a locally finitely presentable Grothendieck category. Then $X \in \G$ is \emph{fp-injective} if $\Ext^1_\G(F,X)=0$ for all $F \in \fp\G$. We denote the full subcategory of $\G$ formed by the fp-injective objects by $\Fpinj\G$.
\end{defn}

\begin{warn} \label{warn:ext}
In order to prevent confusion, we stress that the $\Ext^n_\G$ in this and the next section stands for the \emph{Ext-group \wrt the abelian exact structure} rather than \wrt the pure exact structure. Similarly, $\Ext^n_{\Cpx(\G)}$ will denote the Ext-groups for $\Cpx(\G)$ \wrt the abelian exact structure.
\end{warn}

For basic properties of the class of fp-injective objects and their relation to the injective objects we refer to Appendix~\ref{sec:fp-inj}.
As mentioned in~\cite{Ste70}, replacing injective objects by fp-injective ones can (and does) prove very useful when trying to generalize homological results from a noetherian to a non-noetherian setting. The advantage of fp-injective objects in our case is that on one hand finitely presentable objects see them as injective, and on the other hand the class of fp-injective objects is closed under coproducts.

As $\Fpinj\G$ carries a natural exact structure inherited from $\G$, we may consider its unbounded derived category $\Der(\Fpinj\G)$. It is an easy observation that acyclic complexes in $\Cpx(\Fpinj\G)$ are precisely the pure acyclic complexes in $\G$ with fp-injective components. We shall therefore denote this class by $\Cpac(\Fpinj\G)$.

Since triangulated categories like $\Der(\Fpinj\G)$ usually do not have other interesting colimits than coproducts, they cannot be locally finitely presentable. However, there is a well known remedy---the concept of a compactly generated triangulated category (see~\cite{Nee01,Kr05} for a thorough discussion and generalizations).

\begin{defn} \label{defn:comp-gen}
An object $C$ of a triangulated category $\T$ with small coproducts is \emph{compact} if for each collection $(X_i \mid i\in I)$ of objects of $\T$, the canonical morphism
\[ \coprod_{i \in I} \T(C, X_i) \la \T\big(C, \coprod_{i \in I} X_i\big) \]
is an isomorphism. The category $\T$ is \emph{compactly generated} if there exists a small set $\clS \subseteq \T$ of compact objects \st for each $0 \ne X \in \T$ there is a morphism $0 \ne f\colon \Sigma^n S \to X$ for some $S \in \clS$ and $n \in \Z$.
\end{defn}

After the lengthy initial discussion, the following proposition becomes rather easy.

\begin{prop} \label{prop:comp-gen}
Let $\G$ be a locally finitely presentable Grothendieck category. Then $\Der(\Fpinj\G)$ is a compactly generated triangulated category.
\end{prop}

\begin{proof}
The proof is based on the following two observations. First, each object $X \in \G$ admits an augmented proper right fp-injective resolution
\[ A_X\colon \quad 0 \to X \to J^0 \to J^1 \to J^2 \to \cdots \]
in the sense of~\cite[Definition 2.24]{Mur}. That is, $J^n \in \Fpinj\G$ for all $n \ge 0$ and $\HOM_\G(A_X,J)$ is acyclic for each $J \in \Fpinj\G$.
Here we use the functorially complete cotorsion pair $(\X,\Fpinj\G)$ in $\G$ which is provided by Proposition~\ref{prop:complete-cotor-A}, and we construct the augmented resolution using the corresponding approximation sequences (see also~\cite[Remark 2.25]{Mur}). The truncated complex
\[ J_X\colon \quad 0 \to J^0 \to J^1 \to J^2 \to J^3 \to \cdots \]
with $J^0$ in degree $0$ is then a so-called proper right fp-injective resolution, and there is an obvious map $j\colon X \to J_X$ in $\Cpx(\G)$ whose cone is $A_X$.

Second, if $F \in \G$ is finitely generated, then a proper right fp-injective resolution $J_F$ of $F$ is compact in $\Htp(\Fpinj\G)$ by~\cite[Meta-Theorem 2.29]{Mur} and the discussion above it. This is since $j^*\colon \Htp(\G)(J_X,J') \to \Htp(\G)(X,J')$ is an isomorphism for all $J' \in \Htp(\Fpinj\G)$, and since $\Htp(\Fpinj\G)$ has coproducts which are computed component wise by Lemma~\ref{lem:fp-inj-closure}.

Now we fix proper right fp-injective resolutions $J_F$, where $F$ is running over all representatives of isoclasses of finitely presentable objects in $\G$. Let $\clS \subseteq \Cpx(\Fpinj\G)$ by the collection of all such $J_F$. This is by the above discussion a set of compact objects on $\Htp(\Fpinj\G)$ such that the class
\[ \clS^\perp = \{ X \in \Htp(\Fpinj\G) \mid \Htp(\Fpinj\G)(J_F,X) = 0 \textrm{ for each } J_F \in \clS \} \]
consists precisely of all $X \in \Htp(\Fpinj\G)$ \st $\HOM_\G(F,X)$ is acyclic for each $F \in \fp\G$. Since $\fp\G$ is a generating set for $\G$, each such $X$ is acyclic. Moreover, each such $X$ is pure acyclic by Lemma~\ref{lem:char-pure-mono-epi}. A moment's thought reveals that we even have equality $\clS^\perp = \Hpac(\Fpinj\G)$.

It is now a standard result that the Verdier quotient $\Htp(\Fpinj\G)/\clS^\perp$, which is none other than $\Der(\Fpinj\G) \simeq \Htp(\Fpinj\G)/\Hpac(\Fpinj\G)$, is compactly generated triangulated. We refer to~\cite[Ch. 9]{Nee01} or \cite[\S5]{Kr10} for details.
\end{proof}

Our next goal is to show that $\Der(\Fpinj\G)$ is canonically triangle equivalent to $\Htp(\Inj\G)$. Following the philosophy of the paper, we shall do so by constructing model structures on $\Cpx(\G)$ for both the triangulated categories. Here we start by proving the existence of two cotorsion pairs which we will use for the model of $\Der(\Fpinj\G)$, finishing the entire argument later in~\S\ref{subsec:equiv}.

\begin{lem} \label{lem:fp-inj-model-cp}
Let $\G$ be a locally finitely presentable Grothendieck category. Then $\Cpx(\G)$ admits functorially complete cotorsion pairs
\[ (\T\Cof,\Cpx(\Fpinj\G)) \quad \textrm{and} \quad (\Cof,\Cpac(\Fpinj\G)). \]
\end{lem}

\begin{proof}
Note that by Lemma~\ref{lem:ext-adj}(1), a complex $X \in \Cpx(\G)$ belongs to $\Cpx(\Fpinj\G)$ \iff $\Ext^1_{\Cpx(\G)}(D^n(F),X) = 0$ for all $n \in \Z$ and $F \in \fp\G$. Hence Proposition~\ref{prop:complete-cotor-A} yields the existence of a functorially complete cotorsion pair $(\T\Cof,\Cpx(\Fpinj\G))$ in $\Cpx(\G)$.

Similarly, note that $X \in \Cpac(\G)$ \iff $\Ext^1_{\Cpx(\G)}(\Sigma^n F,X) = 0$ for all $n \in \Z$ and $F \in \fp\G$. Indeed, since $D^n(F)$ is always an extension of $\Sigma^{-n} F$ and $\Sigma^{-n-1} F$, a complex which satisfies $\Ext^1_{\Cpx(\G)}(\Sigma^n F,X) = 0$ for all $n$ and $F$ has all components fp-injective. Thus,
\[ 0 = \Ext^1_{\Cpx(\G)}(\Sigma^n F,X) \cong \Ext^1_{\Cpx(\G),\dw}(\Sigma^n F,X) \cong \Htp(\G)(\Sigma^{n-1} F,X), \]
or in other words $\HOM_\G(F,X) = 0$ for all $F \in \fp\G$, and so $X \in \Cpac(\Fpinj\G)$. If we conversely start with $X \in \Cpac(\Fpinj\G)$, then $\Ext^1_{\Cpx(\G)}(\Sigma^n F,X) = 0$ for all $n$ and $F$ follows simply by retracing the steps. Hence the functorially complete cotorsion pair $(\Cof,\Cpac(\Fpinj\G))$ in $\Cpx(\G)$ again exists by Proposition~\ref{prop:complete-cotor-A}.
\end{proof}

\begin{rem} \label{rem:fp-inj-model-coh}
If $\G$ is locally coherent, the two cotorsion pairs from Lemma~\ref{lem:fp-inj-model-cp} are necessarily hereditary and their mere existence implies the existence of a model structure on $\Cpx(\G)$ whose homotopy category is $\Der(\Fpinj\G)$; see \cite{Gil14}. However, it seems far from being obvious from there that the class of trivial objects is the same as for the Positselski-Becker coderived model structure from \cite[Proposition 1.3.6]{Beck14}. Proving this equality is one of the highlights of Section~\ref{sec:coderived} (see Theorem~\ref{thm:fp-inj-model}).
\end{rem}

\subsection{Coacyclic complexes}
\label{subsec:coder}

As the usual derived category $\Der(\G)$ is obtained from $\Cpx(\G)$ by killing acyclic complexes, $\Htp(\Inj\G)$ can be obtained by killing the smaller class of so-called coacyclic complexes; \cite{Pos11}. Let us provide a few details on this process.

Analogously to Baer's criterion for injectivity for module categories, $\G$ admits a set $\clS$ of finitely generated objects such that
\[ \Inj\G = \{ Y \in \G \mid \Ext^1_\G(S,Y) = 0 \textrm{ for each } S \in \clS \}. \]
It is a direct consequence of Lemma~\ref{lem:ext-adj} that
\[ \Cpx(\Inj\G) = \{ Y \in \Cpx(\G) \mid \Ext^1_{\Cpx(\G)}(D^n(S),Y) = 0 \textrm{ for each } S \in \clS \textrm{ and } n \in \Z \}. \]
Hence $\Cpx(\Inj\G)$ is the right hand side of a functorially complete cotorsion pair in $\Cpx(\G)$ by Proposition~\ref{prop:complete-cotor-A}. The left hand side of this cotorsion pair is the subclass of all acyclic complexes which we are looking for.

\begin{defn}[\cite{Pos11,Beck14}] \label{defn:coacyclic}
A complex $X \in \Cpx(\G)$ is called \emph{coacyclic}, if $\Ext^1_{\Cpx(\G)}(X,Y) = 0$ for each $Y \in \Cpx(\Inj\G)$. We denote the class of all coacyclic complexes by $\Ccoac(\G)$.
\end{defn}

\begin{rem} \label{rem:coacyclic}
This definition perfectly agrees with the terminology of \cite{Beck14}. For a locally noetherian $\G$ it matches also the original and somewhat different definition of Positselski; see~\cite[\S3.7]{Pos11}.
\end{rem}

Invoking a result of Becker and Positselski, we have the following.

\begin{prop} \label{prop:coderived}
Let $\G$ be a locally finitely presentable Grothendieck category. Then there is an abelian hereditary model structure on $\Cpx(\G)$ \st
\begin{enumerate}
 \item Every complex is cofibrant.
 \item The trivial complexes are precisely the coacyclic ones.
 \item The fibrant complexes are the ones with injective components.
\end{enumerate}
Moreover, $\Ho\Cpx(\G)$ is triangle equivalent to $\Htp(\Inj\G)$.
\end{prop}

\begin{proof}
The existence of the model structure a special case of \cite[Proposition 1.3.6]{Beck14} while the last part follows from \cite[Proposition 1.1.13]{Beck14}. In fact, it is not difficult to prove the proposition directly from Proposition~\ref{prop:exact-model}, just taking into account that $(\Ccoac(\G), \Cpx(\Inj\G))$ is a functorially complete cotorsion pair in $\Cpx(\G)$ and $\Ccoac(\G) \cap \Cpx(\Inj\G)$ consists precisely of contractible complexes with injective components.
\end{proof}

Following \cite{Pos11,Beck14}, we can define the \emph{coderived category} $\Dco(\G)$ of $\G$ as the homotopy category of the latter model structure. That is, $\Dco(\G) = \Cpx(\G)[\we\inv]$, where $\we$ consists of all cochain complex homomorphisms whose mapping cones are coacyclic. Then $\Dco(\G) \simeq \Htp(\Inj\G)$ and we can use this insight to construct and understand a model structure for the derived category of fp-injectives.

\subsection{The equivalences}
\label{subsec:equiv}

Here we have collected enough information to proceed to the main results. We start with a general consequence of Theorem~\ref{thm:cotorsion-pac}.

\begin{lem} \label{lem:fpinj-coac-easy}
Let $\G$ be a locally finitely presentable Grothendieck category. Then $\Cpac(\Fpinj\G) \subseteq \Ccoac(\G) \cap \Cpx(\Fpinj\G)$.
\end{lem}

\begin{proof}
We only need to prove that $\Cpac(\Fpinj\G) \subseteq \Ccoac(\G)$. In other words, any given exact sequence $\ep\colon 0 \to I \to E \to X \to 0$, where $I \in \Cpx(\Inj\G)$ and $X \in \Cpac(\Fpinj\G)$, must split. Since $\ep$ must be component wise split and certainly $I \in \Cpx(\PInj\G)$, this follows from Theorem~\ref{thm:cotorsion-pac}.
\end{proof}

This is how far one can get for general locally finitely presentable Grothen\-dieck categories. The crux of comparing $\Htp(\Inj\G)$ with $\Der(\Fpinj\G)$ is that we actually obtain equality between $\Cpac(\Fpinj\G)$ and $\Ccoac(\G) \cap \Cpx(\Fpinj\G)$ if $\G$ is locally coherent. This additional assumption cannot be omitted as will be shown in Example~\ref{expl:non-coh-fail}.

\begin{prop} \label{prop:fpinj-coac-coh}
Let $\G$ be a locally coherent Grothendieck category. Then $\Cpac(\Fpinj\G) = \Ccoac(\G) \cap \Cpx(\Fpinj\G)$.
\end{prop}

\begin{proof}
Let $\clS \subseteq \Fpinj\G$ be a small set of objects such that
\[ \Inj\G = \{ Y \in \Fpinj\G \mid \Ext^1_{\Fpinj\G}(S,Y) = 0 \textrm{ for each } S \in \clS \}. \]
The existence of $\clS$ is proved in Proposition~\ref{prop:fp-inj-baer-coh} and that is where the local coherence of $\G$ is needed. By Lemma~\ref{lem:ext-adj} we also have
\[ \Cpx(\Inj\G) = \{ Y \in \Cpx(\Fpinj\G) \mid \Ext^1_{\Cpx(\Fpinj\G)}(D^n(S),Y) = 0 \; \forall S \in \clS, n \in \Z \}. \]
Here we consider the obvious exact structure on $\Cpx(\Fpinj\G)$ obtained by restriction from either the abelian or the induced pure exact structure on $\Cpx(\G)$. Thus, $\Cpx(\Fpinj\G)$ admits a functorially complete cotorsion pair $(\X,\Cpx(\Inj\G))$ \st $\X$ consists precisely of retracts of $\hat\clS$-filtered objects, where
\[ \hat\clS = \{ D^n(S) \mid S \in \clS \textrm{ and } n \in \Z \}. \]
This follows from Proposition~\ref{prop:complete-cotor-A}.

On the other hand, $(\X,\Cpx(\Inj\G))$ is just a restriction of the cotorsion pair $(\Ccoac(\G),\Cpx(\Inj\G))$ from $\Cpx(\G)$ to $\Cpx(\Fpinj\G)$, so that $\X = \Ccoac(\G) \cap \Cpx(\Fpinj\G)$. Thus our task is to prove that each $X \in \X$ is pure acyclic. Since retracts of pure acyclic complexes are pure acyclic, we can assume by the above discussion that $X$ has an $\hat\clS$-filtration $(X_\alpha \mid \alpha \le \lambda)$ (we view the maps in the filtration as inclusions).

Now it is easy to prove that each $X_\alpha$ is pure acyclic by transfinite induction on $\alpha$. For $\alpha = 0$ this is trivial and for limit ordinals this follows from the fact that direct limits of pure exact sequences are pure exact. If $\alpha = \beta+1$ is an ordinal successor, then we have a component wise pure exact sequence of complexes
\[ 0 \la X_\beta \la X_\alpha \la D^{n_\alpha}(S) \la 0 \]
with $S \in \clS$. Since its both ends are pure acyclic ($X_\beta$ by inductive hypothesis), so is $X_\alpha$. The conclusion follows by considering $X = X_\lambda$.
\end{proof}

Now we can finish the discussion of the existence of a nice model for $\Der(\Fpinj\G)$ which was initiated in~\S\ref{subsec:fpinj}. We also confirm \cite[Conjecture~1]{Kr12} not only for coherent rings, but also for locally coherent Grothendieck categories.

\begin{thm} \label{thm:fp-inj-model}
Let $\G$ be a locally coherent Grothendieck category. Then there is an abelian hereditary model structure on $\Cpx(\G)$ represented by the triple of classes $(\Cof,\We,\Fib)$, where
\begin{enumerate}
 \item $\Fib = \Cpx(\Fpinj\G)$, the class of all complexes with fp-injective components,
 \item $\Fib \cap \We = \Cpac(\Fpinj\G)$, the class of pure acyclic complexes of fp-injectives,
 \item $\We = \Ccoac(\G)$, the class of coacyclic complexes.
\end{enumerate}
In particular, we have triangle equivalences
\[ \Ho\Cpx(\G) = \Dco(\G) \simeq \Der(\Fpinj\G) \simeq \Htp(\Inj\G). \]
\end{thm}

\begin{proof}
Consider the two cotorsion pairs 
\[ (\T\Cof,\Cpx(\Fpinj\G)) \quad \textrm{and} \quad (\Cof,\Cpac(\Fpinj\G)). \]
from Lemma~\ref{lem:fp-inj-model-cp}. Since $\G$ is locally coherent, their right hand sides are closed under taking cosyzygies and so they are hereditary. In order to prove that there is a (unique) hereditary abelian model structure satisfying (1)--(3), we are only required to check that
\[\T\Cof = \Cof \cap \Ccoac(\G) \quad\textrm{and}\quad \Cpac(\Fpinj\G) = \Ccoac(\G) \cap \Cpx(\Fpinj\G). \]
The second equality has been proved in Proposition~\ref{prop:fpinj-coac-coh}. To show the first one, note that clearly $\T\Cof \subseteq \C$ and $\T\Cof \subseteq \Ccoac(\G)$ (in terms of the right hand sides of the corresponding cotorsion pairs, this just says that $\Cpx(\Fpinj\G) \supseteq \Cpac(\Fpinj\G)$ and $\Cpx(\Fpinj\G) \supseteq \Cpx(\Inj\G)$). Suppose conversely that $X \in \Cof \cap \Ccoac(\G)$ and consider an approximation sequence
\[ \ep\colon \quad 0 \la J \la T \la X \la 0 \]
with $J \in \Cpx(\Fpinj\G)$ and $T \in \T\Cof$. Since both $T,X \in \Ccoac(\G)$, we have $J \in \Ccoac(\G) \cap \Cpx(\Fpinj\G) = \Cpac(\Fpinj\G)$. As then $\Ext^1_{\Cpx(\G)}(X,J) = 0$, $\ep$~splits and $X \in \T\Cof$. Thus, an abelian hereditary model structure with the required properties exists.

As $\Ho\Cpx(\G)$ depends only on the class of trivial objects, $\Ho\Cpx(\G) = \Dco(\G)$ is just the definition of the coderived category. The equivalence $\Dco(\G) \simeq \Htp(\Inj\G)$ comes from Proposition~\ref{prop:coderived}. Finally, it follows from \cite[Proposition 5.2 and Corollary 5.4]{Gil11} that we can restrict the model structure from the present theorem to $\Cpx(\Fpinj\G)$ and the homotopy category remains canonically equivalent. As a result of the restriction, we obtain a hereditary exact model structure on $\Cpx(\Fpinj\G)$ whose class of trivial objects is precisely $\Ccoac(\G) \cap \Cpx(\Fpinj\G) = \Cpac(\Fpinj\G)$. That is, the weak equivalences are precisely pure quasi-isomorphisms and the homotopy category is $\Der(\Fpinj\G)$.
\end{proof}

In more algebraic terms, what the theorem says can be restated in complete analogy with \cite[Theorem 1.1(2)]{Kr05} and vastly improves \cite[Theorem 3.13]{Nee14} for locally coherent Grothendieck categories.

\begin{cor} \label{cor:fp-inj}
Let $\G$ be a locally coherent Grothendieck category. Then $\Htp(\Inj\G)$ is a compactly generated triangulated category. Moreover, a functor
\[ I\colon \Db(\fp\G) \la \Htp(\Inj\G) \]
which assigns to each bounded complex $X \in \Db(\fp\G)$ a fixed injective resolution $I(X)$ of $X$ and similarly for morphisms, restricts to an equivalence between $\Db(\fp\G)$ and the full subcategory $\Htp(\Inj\G)^c$ of compact objects in $\Htp(\Inj\G)$.
\end{cor}

\begin{proof}
If we combine the triangle equivalence $\Der(\Fpinj\G) \simeq \Htp(\Inj\G)$ with Proposition~\ref{prop:comp-gen}, we immediately get that the (representatives of isoclasses of) injective resolutions of finitely presentable objects of $\G$ form a compact generating set for $\Htp(\Inj\G)$. Since the functor $I$ is well known to be triangulated and fully faithful, and since $\fp\G$ generates $\Db(\fp\G)$ as a triangulated category, the conclusion follows from~\cite[Lemma 2.2]{Nee92}.
\end{proof}

\subsection{The non-coherent case}
\label{subsec:non-coh}

Let us also briefly discuss what happens if $\G$ is not locally coherent. The main purpose is to illustrate how badly the results in~\S\ref{subsec:equiv} may fail. What we still have is the following.

\begin{lem} \label{lem:loc}
Let $\G$ be a locally finitely presentable Grothendieck category. Then $\Htp(\Inj\G)$ is equivalent to $\Der(\Fpinj\G)/\clL$ for a suitable localizing subcategory (i.e.\ a triangulated subcategory closed under coproducts) $\clL \subseteq \Der(\Fpinj\G)$.
\end{lem}

\begin{proof}
By a straightforward variation of results from~\cite[\S5]{Gil11}, one can restrict the model structure from Proposition~\ref{prop:coderived} to $\Cpx(\Fpinj\G)$ without changing the homotopy category. It follows that $\Htp(\Inj\G) \simeq \Dco(\G)$ is obtained from $\Cpx(\Fpinj\G)$ by inverting all morphisms with mappings cones in $\Ccoac(\G) \cap \Cpx(\Fpinj\G)$, whereas $\Der(\Fpinj\G)$ is obtained from $\Cpx(\Fpinj\G)$ by inverting the morphisms with mapping cones in $\Cpac(\Fpinj\G)$. The rest follows from Lemma~\ref{lem:fpinj-coac-easy} and, explicitly, $\clL$ can be taken as the class of all coacyclic complexes in $\Der(\Fpinj\G)$.
\end{proof}

However, the following example shows that the above localization functor $Q\colon \Der(\Fpinj\G) \to \Htp(\Inj\G)$ may not be an equivalence, and in particular that the first part of \cite[Conjecture 1]{Kr12} fails as stated. It also shows that we need to assume local coherence in Proposition~\ref{prop:fpinj-coac-coh}.

\begin{expl} \label{expl:non-coh-fail}
Let $k$ be a field, $\m$ be an infinite dimensional $k$-vector space, and let $R = k \ltimes \m$ be the trivial extension of $k$ by $\m$. One can quickly verify that $R$ is a perfect non-coherent local ring, where $\m$ is the unique maximal ideal.

Suppose that $X$ is a non-projective $R$-module and $P_1 \overset{g}\to P_0 \to X \to 0$ is its minimal projective resolution. Then $\Img g$ is contained in the Jacobson radical of $P_0$, hence it is a semisimple $R$-module. That is, $\Img{g} \cong (R/\m)^{(I)}$ for non-empty set $I$ and consequently $\Ker g \cong \m^{(I)} \cong (R/\m)^{(J)}$. Hence the minimal second syzygy of $X$ is always an \emph{infinitely} generated semisimple $R$-module. Given any semisimple $R$-module $S$, one also quickly checks that
\[ \Ext^2_R(X,S) \cong \Hom_R(\Ker g,S) \cong \Hom_R((R/\m)^{(J)},S) \]
If we substitute $S = \Ker g \cong (R/\m)^{(J)}$, we immediately see that $\Ext^2(X,-)$ does not commute with coproducts since the image of $1_{\Ker g} \colon \Ker g \to \Ker g$ cannot be contained a finite subcoproduct of $(R/\m)^{(J)}$.

It follows from~\cite[Meta-Theorem 2.29]{Mur} that an injective resolution $I(X)$ of $X$ cannot be compact in $\Htp(\Inj\G)$. Indeed, if $I(S)$ is an injective resolution of $S$, then $I(S)$ is a coproduct of $J$ copies of $I(R/\m)$ in $\Htp(\InjR)$ and we have
\[ \Ext^2_R(X,S) \cong \Htp(\G)(X,\Sigma^2 I(S)) \cong \Htp(\G)(I(X),\Sigma^2 I(S)). \]
Hence the canonical localization functor $\Der(\FpinjR) \to \Htp(\InjR)$, which sends proper right fp-injective resolutions of non-projective finitely presentable modules to injective resolutions of the same modules, cannot be an equivalence since it sends compact objects to non-compact ones.

Pushing this idea a little further, we can demonstrate a complex which is in $\Ccoac(\G) \cap \Cpx(\Fpinj\G)$ but not in $\Cpac(\Fpinj\G)$. Namely, consider a quasi-isomorphism 
\[ f\colon I(R/\m)^{(J)} \to I\big((R/\m)^{(J)}\big).\]
Its mapping cone is acyclic and bounded below, hence coacyclic. Note that the domain of $f$ is a coproduct of $J$ copies of $I(R/\m)$ in $\Der(\FpinjR)$ and the codomain is a coproduct of the same $J$ copies in $\Htp(\InjR)$. If $X$ is non-projective finitely presentable, then $\Htp(\G)(X,\Sigma^2 f)$ cannot be an isomorphism by the above discussion, so that the mapping cone of $f$ is not pure acyclic.
\end{expl}

\section{Singularity categories and Krause's recollement}
\label{sec:singularity}

We conclude with a discussion of the so-called singularity category of a locally coherent Grothendieck category (Definition~\ref{defn:sing}). We do not strive for the utmost generality but rather for a reasonably simple criterion under which the recollement situation from~\cite{Kr05} still arises.

First of all note that there is a canonical localization functor $Q\colon \Dco(\G) \to \Der(\G)$. The coderived category is constructed from $\Cpx(\G)$ by formally inverting the morphisms whose mapping cones are coacyclic, while $\Der(\G)$ by inverting all quasi-isomorphisms. This $Q$ is well known to have a fully faithful right adjoint $R$. From the model theoretic point of view, we can denote by $\Cpx(\G)^\mathrm{co}$ the category of complexes with the model structure from Proposition~\ref{prop:coderived} and by $\Cpx(\G)^\mathrm{der}$ the same category of complexes, but with the standard injective model structure for the derived category (see e.g.~\cite[\S2.3]{Hov99} or~\cite[Example 3.2]{Hov02-models-cot} and references there). Then the identity functor $1_{\Cpx(\G)}\colon \Cpx(\G)^\mathrm{der} \to \Cpx(\G)^\mathrm{co}$ is right Quillen and $R$ is its right derived functor.

Alternatively from the algebraic point of view, we start with the functorially complete cotorsion pair $(\Cac(\G),\dgInj\G)$ in $\Cpx(\G)$. Such a pair certainly exists by \cite[Theorem 4.12]{Gil07}, or by \cite[Theorem 4.2(2)]{St13} together with Proposition~\ref{prop:complete-cotor-A}. Then $\Dco(\G) \simeq \Htp(\Inj\G)$ by Proposition~\ref{prop:coderived} and $\Der(\G)$ is equivalent to the full subcategory $\dgInj\G \subseteq \Htp(\Inj\G)$. The subcategory $\dgInj\G$ is reflective in $\Htp(\Inj\G)$ by \cite{BEIJR12} and the reflection corresponds to the localization functor $Q$ under the identifications in place. The right adjoint $R$ is then simply identified with the inclusion $\dgInj\G \subseteq \Htp(\Inj\G)$.

More interestingly, under appropriate conditions we also have a left adjoint $L$ to $Q$. This has been observed in \cite{Kr05} for certain locally noetherian Grothendieck categories and in \cite[Corollary 2.2.2]{Beck14} in the context of dg modules. However, the left adjoint does not exists unconditionally as was demonstrated in~\cite[Example 4.1]{Nee14}. Following a suggestion of Gillespie from~\cite{Gil14-gor-inj}, the theory developed by Krause \cite{Kr05} turns out to work under the following hypothesis.

\begin{hyp} \label{hyp:recollement}
$\G$ is a locally coherent Grothendieck category which admits a small generating set $\clS \subseteq \fp\G$ of objects of finite projective dimension. That is, for each $S \in \clS$ there exists $i_S \gg 0$ such that $\Ext^{i_S}_\G(S,-) \equiv 0$.
\end{hyp}

This hypothesis is satisfied under various circumstances. The simplest case is when $\G$ is a module category. It also holds for the category $\G$ of quasi-coherent sheaves over a quasi-projective scheme with coherent structure sheaf; see~\cite[Examples 2.1.2(c) and Proposition B.8]{Th90}. Finally, such a hypothesis is a must for Gorenstein categories---compare to~\cite[Definition 2.18]{EEGR08}.

The following is an observation due to Gillespie~\cite{Gil14-gor-inj}, analogous to the argument for~\cite[Lemma 2.2]{Kr05} (compare also to~\cite[Corollary 1.2]{Kr05}---there is no mystery there anymore).

\begin{lem} \label{lem:inj-ac}
Let $\G$ and $\clS$ be as in Hypothesis~\ref{hyp:recollement}, and let $X \in \Cpx(\Inj\G)$. Then the following are equivalent:
\begin{enumerate}
\item $X$ is acyclic.
\item $\Htp(\G)(\Sigma^nS,X) = 0$ for all $S \in \clS$ and $n \in \Z$.
\end{enumerate}
\end{lem}

\begin{proof}
Since $\clS$ is a generating set, each complex satisfying (2) is acyclic. On the other hand, suppose that $X$ is acyclic and $S \in \clS$ and $i_S > 0$ are \st $\Ext^{i_S}_\G(S,-) \equiv 0$. Then $\Htp(\G)(S,X) \cong \Ext^{i_S}_\G(S,Z^{-i_S}(X)) = 0$ and similarly for (de)suspensions of $\clS$.
\end{proof}

As an immediate corollary, we obtain another functorially complete cotorsion pair (compare to \cite[Remark 2.2.5]{Beck14}).

\begin{cor} \label{cor:inj-ac}
Let $\G$ be as in Hypothesis~\ref{hyp:recollement}. Then there is a functorially complete cotorsion pair $(\We,\Cac(\Inj\G))$ in $\Cpx(\G)$.
\end{cor}

\begin{proof}
Let $\U \subseteq \G$ be a small set \st
\[ \Cpx(\Inj\G) = \{ Y \in \Cpx(\G) \mid \Ext^1_{\Cpx(\G)}(D^n(U),Y) = 0 \textrm{ for each } U \in \U \textrm{ and } n \in \Z \}. \]
(see \S\ref{subsec:coder}). Given $\clS$ from Hypothesis~\ref{hyp:recollement}, the the $\Ext^1_\G$-orthogonal class to
\[ \{ D^n(U) \mid U \in \U, n \in \Z \} \cup \{ \Sigma^n S \mid S \in \clS, n \in \Z \} \]
equals $\Cac(\Inj\G)$ by the former lemma. Hence there is a functorially complete cotorsion pair as in the statement by Proposition~\ref{prop:complete-cotor-A}.
\end{proof}

More importantly, the hypothesis ensures that the derived category of $\Der(\G)$ is compactly generated by $\clS$.

\begin{prop} \label{prop:comp-gen-D(G)}
Let $\G$ and $\clS$ be as in Hypothesis~\ref{hyp:recollement}. Then $\clS$ is a set of compact generators for $\Der(\G)$.
\end{prop}

\begin{proof}
The fact that $\clS$ is a set of generators is straightforward. Any $X \in \Der(\G)$ can be fibrantly replaced by a complex $Y \in \dgInj(\G)$. Since we have $\Der(\G)(\Sigma^n S,X) \cong \Htp(\G)(\Sigma^n S,Y)$, it follows that $\Der(\G)(\Sigma^n S,X) = 0$ for all $S \in \clS$ and $n \in \Z$ \iff $X$ is acyclic.

A little trickier is to prove that an object $S \in \clS$ is compact in $\Der(\G)$. To this end, fix $i>0$ \st $\Ext^i_\G(S,-) \equiv 0$ and consider the standard $t$-structure $(\Dle0,\Dge0)$ on $\Der(\G)$ (see~\cite[Exemples 1.3.2]{BBD82}). We first claim that $\Der(\G)(S,\Dle{-i}) = 0$. Indeed, let $X \in \dgInj\G$ represent an object of $\Dle{-i}$. Then again
\[ \Der(\G)(S,X) \cong \Htp(\G)(S,X) \cong \Ext^i_{\G}(S,Z^{-i}(X)) = 0, \]
proving the claim.

Next suppose that we have a map $f\colon S \to \coprod_{j \in J} X_j$ in $\Der(\G)$. Since the $t$-structure truncations commute with coproducts, the previous claim allows us replace all $X_j$ by $\tau^{>-i}X_j$. Hence we can without loss of generality assume that $X_j \in \Dge{1-i}$ for all $j$, i.e.\ that all $X_j$ have non-zero components only in degrees $\ge 1-i$.

The last step of the proof consists in showing that the right adjoint $R\colon \Der(\G) \to \Dco(\G)$ preserves the coproduct $\coprod_{j \in J} X_j$. To see that, note that the fibrant replacement of $\coprod_{j \in J} X_j$ in $\Cpx(\G)^\mathrm{der}$ can be constructed classically using~\cite[Lemma 4.6(i)]{Har66}. That is, we have a short exact sequence
\begin{equation} \label{eq:fib-res}
0 \la \coprod_{j \in J} X_j \la I \la A \la 0, 
\end{equation}
where $I$ is a bounded below complex injective components and $A$ is bounded below and acyclic. Now \cite[Theorem 8.4.4]{Hir03} and its proof dictate that $R(\coprod_{j \in J} X_j) \cong I$. Observe also that the bounded below acyclic complex $A$ is also coacyclic since $0 = \Htp(\G)(A,\Sigma J) \cong \Ext^1_{\Cpx(\G)}(A,J)$ for each $J \in \Cpx(\Inj\G)$. Hence \eqref{eq:fib-res} is also a fibrant replacement of $\coprod_{j \in J} X_j$ in $\Dco(\G)$ and so $R(\coprod_{j \in J} X_j) \cong \coprod_{j \in J} R(X_j)$ canonically.

Once we have this, then, viewing $X$ as an object of $\Dco(\G)$, we have
\[ \Der(\G)\Big(Q(X), \coprod_{j \in J} X_j\Big) \cong \Dco(\G)\Big(X, \coprod_{j \in J} R(X_j)\Big). \] 
Since $X$ is compact in $\Dco(\G)$ by Theorem~\ref{thm:fp-inj-model} and Corollary~\ref{cor:fp-inj}, it follows that $f$ factors through a finite subcoproduct of $\coprod_{j \in J} X_j$, finishing the proof that $X$ is compact in $\Der(\G)$.
\end{proof}

At this point we can use the existing results in \cite{Beck14} to quickly obtain a generalization of~\cite[Theorem 1.1]{Kr05}. First of all, however, we need give a precise definition of a singularity category in our context. We refer to \cite{Bu86} and the introduction of~\cite{Kr05} for the geometric motivation of the terminology and for the representation theoretic significance of the concept.

\begin{defn} \label{defn:sing}
Let $\G$ be a Grothendieck category as in Hypothesis~\ref{hyp:recollement}, and denote by $\Dperf(\fp\G)$ the full subcategory of compact objects of $\Der(\G)$. The notation comes from the fact that in various contexts complexes compact in $\Der(\G)$ are called \emph{perfect}. Since each perfect complex has finitely presentable cohomology by the above proposition, we can view $\Dperf(\fp\G)$ as a thick subcategory of $\Db(\fp\G)$. Then the \emph{singularity category} of $\G$ is defined as the Verdier quotient
\[ \Dsing(\fp\G) = \Db(\fp\G)/\Dperf(\fp\G). \]
\end{defn}

The main motive of \cite{Kr05} was the description of an ``inductive completion'' of the quotient $\Dco(\G)/\Der(\G)$. Following the spirit of this text, we first give a model theoretic background for the situation.

\begin{prop} \label{prop:model-sing}
Let $\G$ be a Grothendieck category satisfying Hypothesis~\ref{hyp:recollement} and let $\We \subseteq \Cpx(\G)$ be the class from Corollary~\ref{cor:inj-ac}. Then there are two abelian model structures on $\Cpx(\G)$ with the class of trivial objects equal to $\We$. The corresponding triples as in Proposition~\ref{prop:exact-model} are
\[ \big(\Cac(\G),\We,\Cpx(\Inj\G)\big) \quad \textrm{and} \quad \big(\Cpx(\G),\We,\Cac(\Inj\G)\big). \]
In particular, we have in both cases a triangle equivalence
\[ \Ho\Cpx(\G) \simeq \Hac(\Inj\G). \]
\end{prop}

\begin{proof}
The first model structure is constructed as the right Bousfield localization $\Cpx(\G)^\mathrm{co}/\Cpx(\G)^\mathrm{der}$ of the model categories $\Cpx(\G)^\mathrm{co}$ and $\Cpx(\G)^\mathrm{der}$ for $\Dco(\G)$ and $\Der(\G)$, respectively. We refer to \cite[Propositions 1.4.2 and 1.5.3]{Beck14} for details. The homotopy category is necessarily equivalent to the Verdier quotient $\Dco(\G)/\Img R$, where $R\colon \Der(\G) \to \Dco(\G)$ is the fully faithful right adjoint functor discussed at the beginning of the section. The second model structure from the statement is constructed in complete analogy with \cite[Proposition 2.2.1]{Beck14}, using the cotorsion pair from Corollary~\ref{cor:inj-ac}. In both cases, the homotopy categories are triangle equivalent to $\Hac(\Inj\G)$ by \cite[Proposition 1.1.14]{Beck14}.
\end{proof}

Now we can state the last main result, generalizing~\cite[Theorem 1.1]{Kr05} and analogous to~\cite[Corollary 2.2.2]{Beck14}.

\begin{thm} \label{thm:singularity}
Let $\G$ be a Grothendieck category satisfying Hypothesis~\ref{hyp:recollement}. Then the localization functor $Q\colon \Dco(\G) \to \Der(\G)$ has both a left adjoint $L$ and a right adjoint $R$. In particular, there is a recollement
\[
\xymatrix@1{
\mathbf{S}(\G) \ar[rr]|(.45)\hole|{i_*} &&
\Dco(\G) \ar@/_3ex/[ll]_{i^*} \ar@/^3ex/[ll]^{i^!} \ar[rr]|(.45)\hole|{Q} &&
\Der(\G) \ar@/_3ex/[ll]_{L} \ar@/^3ex/[ll]^{R}
},
\]
where $\mathbf{S}(\G) = \{ X \in \Dco(\G) \mid Q(X) = 0 \}$ is the full subcategory of $\Dco(\G)$ formed by acyclic complexes. In particular,
\begin{enumerate}
 \item $\mathbf{S}(\G)$ is a compactly generated triangulated category,
 \item $\mathbf{S}(\G)$ is triangle equivalent to $\Hac(\Inj\G)$, and
 \item the full subcategory of compact objects of $\mathbf{S}(\G)$ is equivalent to the idempotent completion of the singularity category $\Dsing(\fp\G)$.
\end{enumerate}
\end{thm}

\begin{proof}
This follows by the same argument as \cite[Corollary 2.2.2]{Beck14}. Indeed, the proof of Proposition~\ref{prop:model-sing} implies that $\mathbf{S}(\G)$ is equivalent to the homotopy category of either of the two abelian model structures
\[ \big(\Cac(\G),\We,\Cpx(\Inj\G)\big) \quad \textrm{and} \quad \big(\Cpx(\G),\We,\Cac(\Inj\G)\big) \]
on $\Cpx(\G)$. Let us denote $\Cpx(\G)$ with these model structures by $\Cpx(\G)^\mathrm{sing}$ and $\Cpx(\G)^\mathrm{i\textrm{-}sing}$, respectively. Then the identity functors $\Cpx(\G)^\mathrm{sing} \to \Cpx(\G)^\mathrm{co}$ and $\Cpx(\G)^\mathrm{i\textrm{-}sing} \to \Cpx(\G)^\mathrm{co}$ are left and right Quillen functors, respectively, which gives the left half of the recollement diagram. The right half is only a formal consequence, see~\cite[Chapter 9]{Nee01} or~\cite[\S4]{Kr10}.
In particular, we obtain a sequence of triangulated functors
\[ \Der(\G) \overset{L}\la \Dco(\G) \overset{i^*}\la \mathbf{S}(\G), \]
which is a so-called localizing sequence. That is, $L$ is fully faithful, $i^*$ is a Verdier localization of $\Dco(\G)$ over $\Img L$, and both $L$ and $i^*$ are left adjoints. Moreover, both $\Der(\G)$ and $\Dco(\G)$ are compactly generated by Theorem~\ref{thm:fp-inj-model}, Corollary~\ref{cor:fp-inj} and Proposition~\ref{prop:comp-gen-D(G)}. Since $L$ has a coproduct preserving right adjoint, $L$ also preserves compactness of objects by~\cite[Lemma 5.4.1]{Kr10}.

Up to equivalence, we are now in the situation of~\cite[Theorem 2.1]{Nee92}. That is, the above sequence restricts to a sequence of functors between the full subcategories of compact objects
\[ \Dperf(\fp\G) \overset{L}\la \Db(\fp\G) \overset{i^*}\la \mathbf{S}(\G)^c, \]
and by~\cite[Theorem 2.1]{Nee92}, the induced functor $\bar{i^*}\colon \Dsing(\fp\G) \to \mathbf{S}(\G)^c$ is fully faithful and each $X \in \mathbf{S}(\G)^c$ is a retract of an object in $\Img \bar{i^*}$. In other words, $\bar{i^*}$ is the idempotent completion functor in the sense of~\cite{BaSch01}. It is also well known that $\mathbf{S}(\G)$ is compactly generated in this situation, see~\cite[Theorem 5.6.1]{Kr10}.
\end{proof}

\begin{rem} \label{rem:singularity}
Upon the identifications of the first two terms in the recollement with subcategories of $\Htp(\Inj\G)$ (see Propositions~\ref{prop:coderived} and~\ref{prop:model-sing}), we obtain a more usual form in which versions of this recollement appear in the literature~\cite{Kr05,Beck14}:
\[
\xymatrix@1{
\Hac(\Inj\G) \ar[rr]|(.45)\hole|{i_*} &&
\Htp(\Inj\G) \ar@/_3ex/[ll]_{i^*} \ar@/^3ex/[ll]^{i^!} \ar[rr]|(.45)\hole|{Q} &&
\Der(\G) \ar@/_3ex/[ll]_{L} \ar@/^3ex/[ll]^{R}
}.
\]
The functor $i_*$ is then simply an inclusion, while $i^*$ and $i^!$ are constructed using the approximation sequences for the cotorsion pairs $(\We,\Cac(\Inj\G))$ (see Corollary~\ref{cor:inj-ac}) and $(\Cac(\Inj\G),\dgInj\G)$, respectively. We refer to \cite{BEIJR12} for details.
This point of view found an interesting application in~\cite{Mur13}.
\end{rem}

We conclude the paper with a short comment on non-noetherian versions of Gorenstein homological algebra and how our techniques apply there. One approach, followed in~\cite{DASS14}, is to consider small additive categories $\clS$ \st $\G = \Modu\clS$ is a locally coherent Gorenstein category in the sense of~\cite{EEGR08}. The results of the present paper will be used in~\cite{DASS14} to show that this leads to a satisfactory generalization of the classical results from~\cite{Bu86}, applies to the theory of triangulated categories, and explains some phenomena in variants of Kasparov's bivariant K-theory of $C^\ast$-algebras.

Another approach was considered by Ding and Chen~\cite{DiCh93,DiCh96} and Gillespie~\cite{Gil10}. They studied two-sided coherent rings with finite self-fp-injective dimension. This is formally a generalization of Iwanaga-Gorenstein rings. As we shall show now, this approach also perfectly fits into our framework. We refer for the background and detailed explanations to~\cite{Gil10}.

\begin{prop} \label{prop:ding-chen}
Let $R$ be a left and right coherent ring of with both ${_RR}$ and $R_R$ of finite fp-injective dimension. Then Gillespie's stable module category, as introduced in~\cite{Gil10}, is triangle equivalent to $\Hac(\InjR)$ and as such it is compactly generated triangulated. 
\end{prop}

\begin{proof}
Upon unraveling~\cite[Definition 3.2]{Gil10} and comparing it with \cite[Theorem 4.7]{Gil10}, we only need to show that a complex of injective modules
\[ I\colon \quad \dots \la I^{-2} \la I^{-1} \la I^0 \la I^1 \la \dots \]
is acyclic \iff $\HOM_R(E,I)$ is acyclic for each fp-injective $R$-module~$E$.

Suppose first that $I$ is acyclic. If $E$ is fp-injective, then the flat dimension of $E$ is finite by~\cite[Theorem 4.2]{Gil10}. Hence
\[ \Htp(\ModR)(E,\Sigma^m I) \cong \Ext^1_R(E,Z^{m-1}(I)) = 0 \]
by Corollary~\ref{cor:gor-inj}, and so $\HOM_R(E,I)$ is acyclic.

Conversely, suppose that $\HOM_R(E,I)$ is acyclic for each fp-injective $R$-module~$E$. By assumption there is an exact sequence $0 \to R \to E_0 \to E_1 \to \cdots \to E_n \to 0$ in $\ModR$ with all $E_0, E_1, \dots, E_n$ fp-injective. Hence there is an exact sequence
\[ 0 \la \HOM_R(E_n,I) \la \cdots \la \HOM_R(E_0,I) \la \HOM_R(R,I) \la 0 \]
of complexes of abelian groups, and, since all $\HOM_R(E_j,I)$ are acyclic for $j = 0, \dots, n$, so is $\HOM_R(R,I)$. The conclusion follows since $\HOM_R(R,I) \cong I$.
\end{proof}

\appendix

\section{Relative filtrations in exact categories}
\label{sec:rel_filtr}

Although exact categories have been studied from various points of view for several decades, this has usually only been done indirectly using the Quillen-Gabriel embedding theorem (see~\cite[Theorem A.1]{Bu10} and the references there). The theorem says that a small exact category is up to equivalence none other than an extension closed subcategory of an abelian category and the conflations are the short exact sequences whose all terms are in the subcategory.

For our purpose, there are two drawbacks of this approach: (a) our exact categories are not small, but more importantly (b) the embedding of an exact category into an abelian one need not preserve infinite colimits. Surprisingly, the choice of literature treating infinite constructions in exact categories other than products or coproducts is very limited. Some steps have been taken in~\cite{Kel90,SaSt11,St14} but that seems to be all. In this appendix, we are going to give a detailed proof for a technical step in Eklof's Lemma (Proposition~\ref{prop:eklof}) which illustrates the difficulty with lack of suitable references.

Let us start with a necessary definition. Morally, $\clS$-filtered objects correspond to absolute cell complexes in homotopy theory (see~\cite[\S1.1]{SaSt11}). What we actually do in the proof of Eklof's Lemma is that we consider the following analog of \emph{relative} cell complexes:

\begin{defn} \label{defn:rel_filtr}
Let $\A$ be an exact category and $\clS$ be a class of objects in~$\A$. A \emph{relative $\clS$-filtration} is a $(\lambda+1)$-sequence $(Y_\alpha, j_{\alpha\beta} \mid \alpha<\beta\le\lambda)$, where $\lambda$ is an ordinal number, \st
\begin{enumerate}
\item $j_{\alpha\beta}\colon Y_\alpha \to Y_\beta$ is an inflation in $\A$ for each $\alpha<\beta\le\lambda$,
\item $\Coker j_{\alpha,\alpha+1} \in \clS$ for each $\alpha<\lambda$.
\end{enumerate}
In other words, the system satisfies conditions (2) and (3) of Definition~\ref{defn:filtr}, but not necessarily condition (1). A morphism $f\colon Z \to X$ in $\A$ is \emph{$\clS$-filtered} if it is a composition of a relative $\clS$-filtration (recall Definition~\ref{defn:lambda-seq}).
\end{defn}

In fact, the latter concept is only auxiliary as the relative filtrations are characterized as follows.

\begin{prop} \label{prop:rel_filtr}
Let $\A$ be an exact category and $\clS$ be a class of objects in~$\A$. The following are equivalent for a morphism $f\colon Z \to X$:
\begin{enumerate}
\item $f$ is $\clS$-filtered (Definition~\ref{defn:rel_filtr}),
\item $f$ is an inflation and $\Coker f$ is $\clS$-filtered (Definition~\ref{defn:filtr}).
\end{enumerate}
\end{prop}

\begin{proof}
The implication $(1) \Rightarrow (2)$ is easy. Starting with a relative $\clS$-filtration $(Y_\alpha, j_{\alpha\beta} \mid \alpha<\beta\le\lambda)$ of $f$, we can construct a direct system $(\Coker j_{0\alpha} \mid \alpha \leq \lambda)$ with the obvious cokernel morphisms. This is clearly an $\clS$-filtration of $\Coker f$. More specifically, it is a $(\lambda+1)$-sequence since direct limits commute with cokernels. Furthermore, condition (1) of Definition~\ref{defn:filtr} is clear, (2) and (3) follow from standard properties of exact categories.

For $(2) \Rightarrow (1)$ there is an easy guess what the relative filtration of $f$ would be: the pullback of a filtration of $\Coker f$ along the deflation $X \to \Coker f$. What does not seem to be so clear is why this pullback is a $(\lambda+1)$-sequence. This is indeed the case and completes the proof of Proposition~\ref{prop:rel_filtr} as we shall see in the lemma below.
\end{proof}

\begin{lem} \label{lem:pull_back_filtr}
Let $\A$ be an exact category, $\clS$ be a class of objects and 
\[\ep\colon \quad 0 \la Y \la E \stackrel{\pi}\la X \la 0\]
be a conflation. Suppose that we are given an $\clS$-filtration $(X_\alpha, i_{\alpha\beta} \mid \alpha<\beta\le\lambda)$ for $X$ and denote for each $\alpha\le\lambda$ by 
\[\ep_\alpha\colon \quad 0 \la Y \stackrel{j_0\alpha}\la E_\alpha \stackrel{\pi_\alpha}\la X_\alpha \la 0 \]
the conflation obtained as a pullback of $\ep$ along the map $i_{\alpha\lambda}$. Then the induced direct system $(E_\alpha, j_{\alpha\beta} \mid \alpha<\beta\le\lambda)$ is a relative $\clS$-filtration of $\pi$.
\end{lem}

\begin{proof}
Note that $E_0 = Y$ and $E_\lambda = E$. Since $\ep_\alpha$ is a pullback of $\ep_\beta$ along $i_{\alpha\beta}$ for each $\alpha<\beta$, conditions~(1) and~(2) of Definition~\ref{defn:rel_filtr} follow immediately from~\cite[Propositions 2.12 and 2.15]{Bu10}.

It remains to prove that $(E_\alpha, j_{\alpha\beta} \mid \alpha<\beta\le\lambda)$ is a $(\lambda+1)$-sequence.
To this end, let $\mu\le\lambda$ be a limit ordinal and let $(f_\alpha\colon E_\alpha \to Z \mid \alpha < \mu)$ be a cocone of $(E_\alpha, j_{\alpha\beta} \mid \alpha<\beta<\mu)$ in $\A$. We must to prove that there is a unique map $f\colon E_\mu \to Z$ satisfying $f\circ j_{\alpha\mu} = f_\alpha$ for each $\alpha<\mu$. The uniqueness of such an $f$ is straightforward. Indeed, suppose that $f,f'\colon E \to Z$ are two such maps, so that $(f-f')\circ j_{\alpha\mu} = 0$ for all $\alpha<\mu$. Specializing to $\alpha=0$, we see that $f-f'$ must factor as $(f-f') = h\circ \pi_\mu$ for some map $h\colon X_\mu \to Z$. But then the equalities
\[ 0 = (f-f') j_{\alpha\mu} = h \pi_\mu j_{\alpha\mu} = h i_{\alpha\mu} \pi_\alpha \]
imply that $h i_{\alpha\mu} = 0$ for all $\alpha<\mu$. Since $X_\mu = \varinjlim_{\alpha<\mu} X_\alpha$, we deduce that $h=0$ and $f=f'$.

Hence we are left with proving the existence of $f$. To this end, we replace the original cocone by $g_\alpha = (^{f_\alpha}_{j_{\alpha\mu}})\colon E_\alpha \to Z \oplus E_\mu$. Clearly it suffices to prove the existence of $g\colon E_\mu \to Z \oplus E_\mu$ \st $g\circ j_{\alpha\mu} = g_\alpha$ for each $\alpha<\mu$. The advantage of the new cocone is that all the maps $g_\alpha$ are inflations by the dual of~\cite[Lemma A.1]{SaSt11}. In particular, the cokernel of $g_0\colon Y \to Z \oplus E_\mu$ exists and for each $\alpha<\mu$ we have a commutative diagram of the form
\begin{equation} \label{eq:seq_eta}
\begin{CD}
\ep_\alpha\colon\quad @. 0 @>>> Y @>{j_{0\alpha}}>> E_\alpha    @>{\pi_\alpha}>> X_\alpha @>>> 0 \\
@.                       @.    @|                 @V{g_\alpha}VV                @VV{t_\alpha}V   \\
\eta\colon\quad       @. 0 @>>> Y @>{g_0}>>      Z \oplus E_\mu @>{\rho}>>           Q    @>>> 0
\end{CD}
\end{equation}
with conflations in rows. It also follows that the maps $(t_\alpha \mid \alpha<\mu)$ defined by the diagram form a cocone of $(X_\alpha, i_{\alpha\beta} \mid \alpha<\beta<\mu)$. Hence we can also construct the colimit map $t\colon X_\mu \to Q$ \st $t\circ i_{\alpha\mu} = t_\alpha$ for each $\alpha<\mu$. Let $p_{E_\mu}\colon Z \oplus E_\mu \to E_\mu$ be the projection onto the second component. Since $p_{E_\mu}\circ g_0 = j_{0\mu}$, we can define a map $q\colon Q \to X_\mu$ using the following commutative diagram
\[
\begin{CD}
\eta\colon\quad       @. 0 @>>> Y @>{g_0}>>      Z \oplus E_\mu @>{\rho}>>           Q    @>>> 0 \\
@.                       @.    @|                 @V{p_{E_\mu}}VV                 @VV{q}V        \\
\ep_\mu\colon\quad    @. 0 @>>> Y @>{j_{0\mu}}>>       E_\mu    @>{\pi_\mu}>>       X_\mu @>>> 0
\end{CD}
\]
For each $\alpha<\mu$ we also have $p_{E_\mu} \circ g_\alpha = j_{\alpha\mu}$ and so also $q \circ t_\alpha = i_{\alpha\mu}$ and, by passing to colimit, $q \circ t = 1_{X_\mu}$. Thus, the pullback of $\eta$ along $t$ equals $\ep_\mu$ and for each $\alpha<\mu$ we have a commutative diagram
\[
\begin{CD}
\ep_\alpha\colon\quad @. 0 @>>> Y @>{j_{0\alpha}}>> E_\alpha    @>{\pi_\alpha}>> X_\alpha @>>> 0 \\
@.                       @.    @|              @V{j_{\alpha\mu}}VV           @VV{i_{\alpha\mu}}V \\
\ep_\mu\colon\quad    @. 0 @>>> Y @>{j_{0\mu}}>>    E_\mu       @>{\pi_\mu}>>    X_\mu    @>>> 0 \\
@.                       @.    @|              @V{\tilde g}VV             @VV{t}V                \\
\eta\colon\quad       @. 0 @>>> Y @>{g_0}>>      Z \oplus E_\mu @>{\rho}>>         Q      @>>> 0 
\end{CD}
\]
Here the map $\tilde g$ in the diagram is defined a the pullback of $t$ along $\rho$.

By comparing the latter diagram with \eqref{eq:seq_eta}, we obtain for each $\alpha<\mu$ a unique map $d_\alpha\colon X_\alpha \to Y$ \st $g_\alpha - \tilde{g} \circ j_{\alpha\mu} = g_0 \circ d_\alpha \circ \pi_\alpha$. It is also straightforward from the construction $d_\beta i_{\alpha\beta} = d_\alpha$ for each $\alpha<\beta<\mu$,. Thus $(d_\alpha \mid \alpha<\mu)$ is a cocone of $(X_\alpha, i_{\alpha\beta} \mid \alpha<\beta<\mu)$ and there is a unique map $d\colon X_\mu \to Y$ \st $d i_{\alpha\mu} = d_\alpha$ for each $\alpha<\mu$. Now we define our map $g\colon E_\mu \to Z \oplus E_\mu$ as $g = \tilde{g} + g_0 d \pi_\mu$. A short computation reveals that
\[ g_\alpha - gj_{\alpha\mu} = g_\alpha - \tilde{g}j_{\alpha\mu} - g_0d\pi_\mu j_{\alpha\mu} = g_0d_\alpha\pi_\alpha - g_0di_{\alpha\mu}\pi_\alpha = 0 \]
for each $\alpha<\mu$. Hence we have found $g\colon E_\mu \to Z \oplus E_\mu$ with the desired property.
\end{proof}

\section{Fp-injective objects in Grothendieck categories}
\label{sec:fp-inj}

A part of this paper relies on properties of fp-injective objects in locally finitely presentable or locally coherent Grothendieck categories, whose analogs for categories of modules over rings (coherent rings, resp.)\ are mostly well known~\cite{Ste70}. As we have not found a suitable reference for our situation, we will discuss it here, using some well known properties of locally finitely presentable Grothendieck categories from \cite{Ste75}. In particular, we will freely use that each object $X \in \G$ has a small set of subobjects and these form a modular and upper-continuous complete lattice by~\cite[Propositions IV.5.3 and V.1.1 (c)]{Ste75}. Recall also that $F \in \G$ is \emph{finitely generated} if $\G(F,-) \colon \G \to \Ab$ preserves direct unions (\cite[\S V.3]{Ste75}), or equivalently if it is a quotient of a finitely presentable object. Another very useful fact to note is that given a short exact sequence
\[ 0 \to K \to L \to M \to 0 \]
in $\G$ with $M$ finitely presentable and $L$ finitely generated, then $K$ is finitely generated; see~\cite[Proposition 3.4]{Ste75}. In particular, $\G$ is locally coherent \iff finitely generated subobjects of finitely presentable objects are always finitely presentable.

Now recall that $X \in \G$ is fp-injective if $\Ext^1_\G(F,X) = 0$ for all $F \in \fp\G$. As in Sections~\ref{sec:coderived} and~\ref{sec:singularity}, $\Ext^n_\G$ denotes here the Ext-group \wrt the abelian exact structure and the class of fp-injective objects is denoted by $\Fpinj\G$. The class $\Fpinj\G$ has certain closure properties, completely analogous to those for module categories.

\begin{lem} \label{lem:fp-inj-closure}
Let $\G$ be a locally finitely presentable Grothendieck category. The class $\Fpinj\G$ is closed under extensions, pure subobjects, direct unions, coproducts and products.
\end{lem}

\begin{proof}
The closure under extensions and products follows directly from the definition. The closure under direct unions and coproducts has been proved in \cite[Corollaries 2.3 and 2.4]{Ste70}. Although the relevant statements in \cite[\S2]{Ste70} are for module categories, their proofs can be taken verbatim for locally finitely presentable Grothendieck categories. Finally, suppose that $0 \to X \to Y \overset{p}\to Z \to 0$ is pure exact in $\G$ and $Y \in \Fpinj\G$. Then we have for each $F \in \fp\G$ an exact sequence
\[ \G(F,Y) \overset{\G(F,p)}\la \G(F,Z) \la \Ext^1_\G(F,X) \la \Ext^1_\G(F,Y) = 0. \]
Since $\G(F,p)$ is surjective, we get $\Ext^1_\G(F,X) = 0$ and $X$ is fp-injective.
\end{proof}

Fp-injective objects are also called \emph{absolutely pure} because of the following property.

\begin{lem} \label{lem:abs-pure}
The following are equivalent for $X \in \G$:
\begin{enumerate}
\item $X$ is fp-injective.
\item Each short exact sequence $0 \to X \overset{i}\to Y \overset{p}\to Z \to 0$ in $\G$ is pure exact.
\end{enumerate}
\end{lem}

\begin{proof}
If $X$ is fp-injective, then $\G(F,p)$ is surjective for each $F \in \fp\G$ since $\Ext^1_\G(F,X) = 0$. Hence (2) holds. Conversely assuming (2), the exact sequence $0 \to X \to E(X) \to C \to 0$, where $E(X)$ is the injective envelope of $X$, is pure. Since clearly $E(X)$ is fp-injective, so is $X$.
\end{proof}

A word of caution is due here. $\Fpinj\G$ is \emph{not} closed under taking direct limits or cosyzygies unless $\G$ is locally coherent. This is also completely analogous to the case of fp-injective modules.

\begin{prop} \label{prop:fp-inj-loc-coh}
Let $\G$ be a locally finitely presentable Grothendieck category. The following are equivalent:

\begin{enumerate}
\item $\G$ is locally coherent.
\item $\Fpinj\G$ is closed in $\G$ under cokernels of monomorphisms.
\item $\Fpinj\G$ is closed in $\G$ under direct limits.
\end{enumerate}
\end{prop}

\begin{proof}
(1) $\implies$ (2). It suffices to prove that $\Ext^2_\G(F,X) = 0$ for each $F \in \fp\G$ and $X \in \Fpinj\G$ since (2) then follows by a simple homological argument. To see that, consider a $2$-extension represented by
\[ \ep\colon \quad 0 \la X \la E_1 \overset{e}\la E_2 \overset{p}\la F \la 0. \]
Using a variation of~\cite[Lemma 13.2.1]{KaSch06} or~\cite[Corollary 5.4]{St14}, we will replace $\ep$ by an equivalent $2$-extension with the $E_2$-term finitely presentable. To this end, the assumption that $\fp\G$ is generating allows us to choose an epimorphism $q\colon \coprod_{i \in I} F_i \to E_2$ with all $F_i$ finitely presentable. In particular, $pq\colon \coprod_{i \in I} F_i \to F$ is an epimorphism and since $F$ is finitely generated, there is a finite subset $I' \subseteq I$ such that the restriction $\coprod_{i \in I'} F_i \to F$ of $pq$ is still an epimorphism. Hence by taking the pullback of $e$ along the restriction $q'\colon \coprod_{i \in I'} F_i \to E_2$ of $q$, we obtain a commutative diagram with exact rows:
\[
\begin{CD}
\ep'\colon \quad 0 @>>> X @>>> E'_1 @>{e'}>> \coprod_{i \in I'} F_i @>{pq'}>> F  @>>>  0   \\
                @.     @|     @VVV                   @VV{q'}V                @|            \\
\ep\colon \quad  0 @>>> X @>>> E_1  @>{e}>>            E_2          @>{p}>>   F  @>>>  0.
\end{CD}
\]
The upper row $\ep'$ is the new $2$-extension which we have been looking for.

Finally, since $\G$ is locally coherent, $\Img e'$ is finitely presentable. Thus, $0 \to X \to E'_1 \to \Img e' \to 0$ splits since $X$ is fp-injective and both $\ep'$ and $\ep$ must represent the zero element of $\Ext^2_\G(F,X)$. Since $\ep$ was chosen arbitrarily, $\Ext^2_\G(F,X) = 0$.

(2) $\implies$ (3). Let $(X_i \mid i \in I)$ be a direct system in $\Fpinj\G$. Then we have a pure exact sequence
\[ 0 \la K \la \coprod_{i \in I} X_i \la \varinjlim_{i \in I} X_i \la 0 \]
in $\G$. Both $K$ and $\coprod_{i \in I} X_i$ are fp-injective by Lemma~\ref{lem:fp-inj-closure}, and so is $\varinjlim_{i \in I} X_i$ by our hypothesis.

(3) $\implies$ (1). The proof of the implication (d)~$\Rightarrow$~(a) of \cite[Theorem 3.2]{Ste70} on page 326 is perfectly valid for locally finitely presentable Grothendieck categories, not only for module categories.
\end{proof}

Now we focus on the properties of $\Fpinj\G$ as an exact category. As usual, conflations are the short exact sequences in $\G$ with all terms fp-injective. In fact, the restrictions of the abelian and the pure exact structures to $\Fpinj\G$ coincide by Lemma~\ref{lem:abs-pure}. We will be first concerned with projective and injective objects in $\Fpinj\G$.

\begin{prop} \label{prop:fp-inj-enough}
Let $\G$ be a locally finitely presentable Grothendieck category and consider $\Fpinj\G$ as an exact category as above.

\begin{enumerate}
\item $\Fpinj\G$ has enough projectives.
\item $\Fpinj\G$ has enough injectives if $\G$ is locally coherent.
\end{enumerate}
\end{prop}

\begin{proof}
(1) By Proposition~\ref{prop:complete-cotor-A}, there is a functorially complete cotorsion pair $(\X,\Fpinj\G)$ in $\G$. Clearly, $\clP = \X \cap \Fpinj\G$ is the class of projective objects in $\Fpinj\G$, and by restriction we get a functorially complete cotorsion pair $(\clP,\Fpinj\G)$ in $\Fpinj\G$.

(2) Given $X \in \Fpinj\G$, consider the usual injective envelope $E(X)$ of $X$. Then $E(X)/X$ is also fp-injective by Proposition~\ref{prop:fp-inj-loc-coh}, and so $X \to E(X)$ is an inflation in $\Fpinj\G$.
\end{proof}

\begin{rem} \label{rem:fp-inj-enough}
It follows that if $\G$ is locally coherent, the injective objects in $\Fpinj\G$ are precisely the injective objects in $\G$. We do not know whether or when $\Fpinj\G$ has enough injective objects if $\G$ is not locally coherent.
\end{rem}

We finish the appendix with a version of Baer's criterion for injectivity for $\Fpinj\G$. The following technical lemma is crucial.

\begin{lem} \label{lem:pure-refine}
Let $\G$ be a locally coherent Grothendieck category. Then there exists a regular cardinal $\kappa$ such that each $X \in \G$ can be filtered by $\kappa$-presentable objects in the pure exact structure on $\G$.
\end{lem}

\begin{proof}
By~\cite[Theorem 2.33]{AR94} and the following remark, there exist arbitrarily large regular cardinals $\kappa$ such that each morphism $f\colon F \to X$ with a $\kappa$-presentable domain $F$ factors as $F \to X' \overset{i}\to X$, where $i$ is a pure monomorphism and $X'$ is $\kappa$-presentable. Note that \cite{AR94} takes the equivalent condition of Lemma~\ref{lem:char-pure-mono-epi}(2) as the definition of pure monomorphisms.

Let us fix one such $\kappa$ with the property above and suppose we have $X \in \G$. We must construct a filtration $(X_\alpha, i_{\alpha\beta} \mid \alpha<\beta\le\lambda)$ of $X$ where all $i_{\alpha\beta}$ are pure monomorphisms and all $\Coker i_{\alpha,\alpha+1}$ are $\kappa$-presentable. In particular we can and will view all the $i_{\alpha\beta}$ as inclusions between subobjects of $X$ in $\G$. We will construct the filtration by induction on $\alpha$, starting according to Definition~\ref{defn:filtr} with $X_0 = 0$. If $\alpha$ is a limit ordinal and all $X_\beta$, $\beta < \alpha$, have been constructed, we must take $X_\alpha = \bigcup_{\beta<\alpha} X_\beta$. Since a direct limit of pure monomorphisms is again a pure monomorphism, the inclusions $X_\beta \subseteq X_\alpha$ for all $\beta < \alpha$ as well as $X_\alpha \subseteq X$ are pure.

Suppose finally that we are at a successor stage. That is, $\alpha = \beta+1$ and $X_\beta$ has been constructed. Then we can assume that $X_\beta$ is a proper subobject of $X$, or else we could finish the construction by putting $\lambda = \beta$. So $X/X_\beta \ne 0$ and there is a non-zero morphism $f\colon F \to X/X_\beta$ with $F \in \fp\G$. By the choice of $\kappa$ we know that $f$ factors through $X' \subseteq X/X_\beta$ where $X'$ is a $\kappa$-presentable pure subobject and necessarily $0 \ne X'$. By pulling back $X' \subseteq X/X_\beta$ along the pure epimorphism $X \to X/X_\beta$, we obtain a pure inclusion $X'' \subseteq X$, and we put $X_\alpha = X''$. Then clearly $X_\alpha/X_\beta \cong X'$ is $\kappa$-presentable and the embeddings $X_\beta \subseteq X_\alpha$ and $X_\alpha \subseteq X$ are pure.
\end{proof}

For our version of Baer's criterion, we will restrict ourselves only to locally coherent Grothendieck categories. The reason is clear from Proposition~\ref{prop:fp-inj-enough}.

\begin{prop} \label{prop:fp-inj-baer-coh}
Let $\G$ be a locally coherent Grothendieck category and consider $\Fpinj\G$ as an exact category. Then there is a small set $\clS \subseteq \Fpinj\G$ of objects \st
\[ \Inj\G = \{ X \in \Fpinj\G \mid \Ext^1_{\Fpinj\G}(S,X) = 0 \textrm{ for each } S \in \clS \}. \]
\end{prop}

\begin{proof}
One can interpret Proposition~\ref{prop:fp-inj-enough} and Remark~\ref{rem:fp-inj-enough} as follows: If $\G$ is locally coherent, there is a functorially complete cotorsion pair $(\Fpinj\G,\Inj\G)$ in $\Fpinj\G$. Indeed, we get such a cotorsion pair by restricting to $\Fpinj\G$ the functorially complete cotorsion pair $(\G,\Inj\G)$ in $\G$.

Now let $\kappa$ be as in Lemma~\ref{lem:pure-refine} and let $\clS$ be a representative set of all $\kappa$-presentable fp-injective objects in $\G$. Consider arbitrary $X \in \Fpinj\G$ and a filtration $(X_\alpha \mid \alpha\le\lambda)$ of $X$ by $\kappa$-presentable objects \wrt the pure exact structure. Then all $X_\alpha$ are fp-injective by Lemma~\ref{lem:fp-inj-closure}, and the factors $X_{\alpha+1}/X_\alpha$ are also fp-injective by Proposition~\ref{prop:fp-inj-loc-coh}. In other words, $\Fpinj\G = \Filt\clS$, so $\Fpinj\G$ is a deconstructible class in itself.

Now $\{ X \in \Fpinj\G \mid \Ext^1_{\Fpinj\G}(S,X) = 0 \textrm{ for each } S \in \clS \} \subseteq \Inj\G$ by Proposition~\ref{prop:eklof}, while the other inclusion is trivial.
\end{proof}

\bibliographystyle{alpha}
\bibliography{derived_fpinj}

\end{document}